\documentclass[aos,MSNbibl,nameyear,seceqn,dvips]{arximspdf}
\usepackage{mathrsfs}
\usepackage{dcolumn}
\usepackage{graphicx}

\doi{10.1214/12-AOS993} 
\volume{40}
\issue{2}
\pubyear{2012}
\firstpage{908}
\lastpage{940}

\makeatletter
\newcolumntype{d}[1]{D{.}{.}{#1}}
\newtheorem{theorem}{Theorem}
\newtheorem{lemma}{Lemma}
\makeatother

\begin{document}
\begin{frontmatter}

\title{Two sample tests for high-dimensional covariance~matrices}
\runtitle{Two sample tests for high-dimensional covariance matrices}

\begin{aug}
\author[A]{\fnms{Jun} \snm{Li}\ead[label=e1]{junli@iastate.edu}}
\and
\author[B]{\fnms{Song Xi} \snm{Chen}\thanksref{aut2}\ead
[label=e2]{csx@gsm.pku.edu.cn}\corref{}}
\runauthor{J. Li and S. X. Chen}
\thankstext{aut2}{Supported by a NSFC Key Grant 111 31002.}
\affiliation{Iowa State University, and Peking
University and Iowa State University}
\address[A]{Department of Statistics\\
Iowa State University\\
Ames, Iowa 50011-1210\\
USA\\
\printead{e1}} 
\address[B]{Guanghua School of Management \\
\quad and Center for Statistical Science \\
Peking University\\
Beijing 100871\\
China\\
and\\
Department of Statistics\\
Iowa State University\\
Ames, Iowa 50011-1210\\
USA\\
\printead{e2}}
\end{aug}

\received{\smonth{7} \syear{2011}}
\revised{\smonth{12} \syear{2011}}

%
\begin{abstract}
We propose two tests for the equality of covariance matrices between
two high-dimensional populations.
One test is on the whole variance--covariance matrices, and the other
is on off-diagonal sub-matrices,
which define the covariance between two nonoverlapping segments of the
high-dimensional random vectors.
The tests are applicable (i) when the data dimension is much larger
than the sample sizes, namely the
``large $p$, small $n$'' situations and (ii)
without assuming parametric distributions for the two
populations. These two aspects surpass the capability of
the conventional likelihood ratio test. The proposed tests can be
used to test on covariances associated with gene ontology terms.
\end{abstract}

%

\end{frontmatter}
%

\section{Introduction}\label{sec1}

Modern statistical data are increasingly high dimensional, but with
relatively small sample sizes. Genetic data typically carry thousands
of dimensions for measurements on the genome. However, due to limited
resources available to replicate study objects, the sample sizes are
usually much smaller than the dimension. 
This is the so-called ``large~$p$, small~$n$'' paradigm. An enduring
interest in Statistics is to know if two populations share the same distribution
or certain key distributional characteristics, for instance the mean or
covariance. The two populations here can refer to two ``treatments'' in
a study. As testing for equality of high-dimensional distributions is
far more challenging than that for the fixed-dimensional data, testing
for equality of key characteristics of the distributions is more
achievable and desirable due to easy interpretation. There has been a
set of research on inference for means of high-dimensional
distributions either in the context of multiple testing,
as in \citet{r43}, \citet{r16},
\citet{r22} and \citet{r26}, or in the context\vadjust{\goodbreak} of simultaneous multivariate
testing as in \citet{r5} 
and \citet{r13}. See also \citet{r28}, \citet{r23} and \citet{r45} for inference on high-dimensional conditional means.

In addition to detecting difference among the population means, there
is a~strong motivation for comparing dependence among components of
random vectors under different treatments, as high data dimensions can
potentially increase the complexity of the dependence. In genomic
studies, genetic measurements, either the micro-array expressions or
the single nucleotide polymorphism (SNP) counts, may have an internal
structure dictated by the genetic networks of living cells. And the
variations and dependence among the measurements of the genes may be
different under different biological conditions and treatments. For
instance, some genes may be tightly correlated
in the normal or less severe conditions, but they can become
decoupled due to certain disease progression; see \citet{r41} for a discussion.

There have been advances on inference for high-dimensional covariance
matrices. The probability limits and the
limiting distributions of extreme eigenvalues of the sample covariance
matrix based on the random matrix theory are
developed in \citet{r2}, \citet{r3}, \citet{r42},
\citet{r29} and \citet{r20},
\citet{r30}, \citet{r6} and others. \citet{r44} and Bickel and Levina (\citeyear{r9}, \citeyear{r10}) proposed
consistent estimators to the population covariance matrices by either
truncation or Cholesky decomposition.
\citet{r21}, \citet{r31} and \citet{r32} considered covariance estimation under factor models.
There are also developments in conducting LASSO-type regularization
estimation of high-dimensional covariances in \citet{r27} and \citet{r38}. Despite these
developments, it is still challenging to transform these results to
test procedures on high-dimensional covariance matrices.

As part of the effort in discovering significant differences between
two high-dimensional distributions, we develop in this paper two-sample
test procedures on high-dimensional covariance matrices.
Let $X_{i1},\ldots,X_{i n_i}$ be an independent and identically
distributed sample drawn
from a $p$-dimensional distribution $F_i$, for \mbox{$i=1$} and 2, respectively.
Here the dimensionality $p$ can be a lot larger than the two sample
sizes $n_1$ and $n_2$ so that $p/n_i \to\infty$.
Let~$\mu_i$ and~$\Sigma_i$ be, respectively, the mean vector and
variance--covariance matrix of the $i$th population.
The primary interest is to test
%
\begin{equation}\label{eq:hypo1}
H_{0 a} \dvtx \Sigma_1 = \Sigma_2 \quad \mbox{versus} \quad  H_{1 a} \dvtx
\Sigma_1
\ne\Sigma_2.
\end{equation}
Testing for the above high-dimensional hypotheses is a nontrivial
statistical problem. Designed for fixed-dimensional data, the conventional
likelihood ratio\vadjust{\goodbreak} test [see \citet{r1} for details] may be used for
the above hypothesis under $p \le\min \{n_1, n_2\}$.
If we let
\[
\bar{X_i}=\frac{1}{n_i}\sum_{j=1}^{n_i}X_{ij}\quad \mbox{and}\quad
Q_i=\sum_{j=1}^{n_i}(X_{ij}-\bar{X_i})(X_{ij}-\bar{X_i})^{\prime
},
\]
then the likelihood ratio (LR) statistic for $H_{0 a}$ is
\[
\lambda_n=\frac{\prod_{i=1}^2|Q_i|^{({1}/{2})n_i}}{|Q|^{({1}/{2})n}}
\frac{n^{({1}/{2})pn}}{\prod_{i=1}^2 n_i^{({1}/{2})p
n_i}},
\]
where $Q=Q_1 +Q_2$ and $n=n_1+n_2$. However, when $p > \min\{n_1, n_2\}
$, at least one of the sample covariance matrices $Q_i/(n_i-1)$ is singular
[\citet{r18}]. This causes the LR statistic $-2 \log(\lambda_n)$
to be
either infinite or undefined, which fundamentally alters the limiting
behavior of the LR
statistic. In an important development, \citet{r4} demonstrated
that even when
$p \le\min\{n_1, n_2\}$ where $\lambda_n$ is properly defined, the
test encounters a power
loss if $p \to\infty$ in such a manner that $p/n_i \to c_i \in(0,1)$
for $i=1$ and $2$. By employing the theory
of large dimensional random matrices, \citet{r4} proposed a
correction to the LR statistic
and demonstrated that the corrected test is valid under $p/n_i \to c_i
\in(0,1)$.
\citet{r39} proposed a test based on a metric that measures the
difference between the two sample covariance
matrices by assuming $p/n_i \to c_i \in[0,\infty)$ and the normal
distributions.
There are also one sample tests for a high-dimensional
variance--covariance $\Sigma$. \citet{r34} and Chen, Zhang and Zhong (\citeyear{r14}) introduced tests
for $\Sigma$ being sphericity and identity for normally distributed
random vectors.
\citet{r35} considered a class of covariance estimators
which are convex sums of $S_n$ and $I_p$ under
moderate dimensionality ($p/n \to c$).
\citet{r11} developed tests for $\Sigma$ having a banded
diagonal structure based on random matrix theory. \citet{r33}
developed a bias-corrected test to examine the significance of the
off-diagonal elements of the residual covariance matrix.
All these tests assume either normality or moderate dimensionality such
that $p/n \to c$ for a finite constant $c$, or both.

We develop in this paper two-sample tests on high-dimensional
variance--covariances without the normality assumption while allowing
the dimension to be much larger than the sample sizes. In addition to
testing for the whole variance--covariance matrices, we propose a test on
the equality of off-diagonal sub-matrices in $\Sigma_1$ and $\Sigma
_2$. The interest
on such a test arises naturally in applications, when we are interested
in knowing
if two segments of the high-dimensional data share the same covariance
between the two treatments.
We will argue in Section \ref{sec3} that the two tests on the whole covariance
and the off-diagonal sub-matrices may be
used collectively to reduce the dimensionality of the testing problem.

This paper is organized as follows. We propose the two-sample test for
the whole covariance matrices in Section \ref{sec2}\vadjust{\goodbreak}
which includes
the asymptotic normality of the test statistic and a power evaluation.
Properties of the test for the off-diagonal sub-matrices are reported
in Section \ref{sec3}.
Results from simulation studies are outlined in Section \ref{sec4}. Section \ref{sec5}
demonstrates how to apply the proposed tests on a gene ontology data
set for acute lymphoblastic leukemia. All technical details are
relegated to Section \ref{sec6}.

\section{Test for high-dimensional variance--covariance}\label{sec2}

The test statistic for the hypothesis (\ref{eq:hypo1}) is formulated
by targeting on $\operatorname{tr}\{(\Sigma_1-\Sigma_2)^2\}$, the
squared Frobenius norm of $\Sigma_1-\Sigma_2$. Although the Frobenius
norm is large in magnitude compared with other matrix norms,
using it for testing brings two advantages. One is that test statistics
based on the norm are relatively easier to be analyzed than those based
on the other norm, which is especially the case when considering the
limiting distribution of the test statistics. The latter renders
formulations of test procedures and power analysis, as we will
demonstrate later. The other advantage is that it can be used to
directly target on certain sections of the covariance matrix as shown
in the next section.
The latter would be hard to accomplish with other norms.

As $\operatorname{tr}\{(\Sigma_1-\Sigma_2)^2\}=\operatorname{tr}(\Sigma_1^2)+\operatorname{tr}(\Sigma
_2^2)-2\operatorname{tr}(\Sigma_1\Sigma_2)$, we will construct estimators for each term.
It is noted that $\operatorname{tr}(S_{n h}^2)$, where $S_{n h}$ is the sample
covariance of the $h$th sample, is a poor estimator of $\operatorname{tr}(\Sigma
_h^2)$ under high dimensionality. The idea is to streamline terms in
$\operatorname{tr}(S_{n h}^2)$ so as to make it unbiased to $\operatorname{tr}(\Sigma_h^2)$ and
easier to analyze in subsequent asymptotic evaluations. We consider
U-statistics of form $\frac{1}{n_h(n_h-1)} \sum_{i \ne j}
(X_{hi}^{\prime}
X_{hj})^2$ which is unbiased if $\mu_h =0$. To account for $\mu_h \ne
0$, we subtract two other U-statistics of order three and four,
respectively, using an approach dated back to Glasser (\citeyear{r24}, \citeyear{r25}).
Specifically, we propose
\begin{eqnarray}\label{eqn2}
A_{n_h}
&=& \frac{1}{n_h(n_h-1)} \sum_{i \ne j} (X_{hi}^{\prime}
X_{hj})^2-\frac{2}{n_h(n_h-1)(n_h-2)} \sum_{i,j,k}^{\star}
X_{hi}^{\prime} X_{hj}X_{hj}^{\prime}
X_{hk}\nonumber\\[-8pt]\\[-8pt]
&&{}+\frac{1}{n_h(n_h-1)(n_h-2)(n_h-3)} \sum_{i,j,k,l}^{\star}
X_{hi}^{\prime} X_{hj}X_{hk}^{\prime} X_{hl}\nonumber
\end{eqnarray}
to estimate $\operatorname{tr}(\Sigma_h^2)$. Throughout this paper we use $\sum
^{\star}$ to denote summation over mutually distinct indices. For
example, $\sum^{\star}_{i,j,k}$ means summation over $\{(i,j,k)\dvtx i\ne
j, j\ne k, k\ne i\}$.
Similarly, the estimator for
$\operatorname{tr}(\Sigma_1 \Sigma_2)$ is
%
\begin{eqnarray}\label{eqn3}
C_{n_1n_2}
&=&\frac{1}{n_1n_2}\sum_{i}\sum_j(X_{1i}^{\prime}X_{2j})^2 -\frac{1}{n_1
n_2 (n_1-1)}\sum_{i, k}^{\star}\sum_{j}X_{1 i}^{\prime}X_{2j}X_{2
j}^{\prime}X_{1 k}\nonumber\\
&&{}- \frac{1}{n_1 n_2 (n_2-1)}\sum_{i, k}^{\star}\sum_{j}X_{2
i}^{\prime}X_{1j}X_{1 j}^{\prime}X_{2 k}\\
&&{} + \frac{1}{n_1 n_2 (n_1-1)(n_2-1)}\sum_{i, k}^{\star}\sum_{j,
l}^{\star}X_{1
i}^{\prime}X_{2j}X_{1 k}^{\prime}X_{2 l}.\nonumber
\end{eqnarray}

There are other ways to attain estimators for $\operatorname{tr}(\Sigma_h^2)$ and
$\operatorname{tr}(\Sigma_1 \Sigma_2)$. In fact, there is a family of estimators for
$\operatorname{tr}(\Sigma_h^2)$ in the form of $\operatorname{tr}(S_h^2)-\alpha_{n_h}\sum
_{i=1}^{n_h}\operatorname{tr}\{(X_{hi}X_{hi}^{\prime}-S_h)^2\}$ where $\alpha_{n_h}
= \alpha/n_h^2$ for any constant $\alpha$. A family can be similarly
formulated for $\operatorname{tr}(\Sigma_1 \Sigma_2)$. It can be shown that this
family of estimators is asymptotically equivalent to the proposed
$A_{n_h}$ in the sense that they share the same leading order term.
However, this family is more complex than the proposed.

The test statistic is
%
\begin{equation}\label{eq:2.2a}
T_{n_1,n_2}=A_{n_1}+A_{n_2} - 2 C_{n_1n_2}
\end{equation}
which is unbiased for $\operatorname{tr}\{(\Sigma_1-\Sigma_2)^2\}$. Besides the
unbiasedness, $T_{n_1,n_2}$ is invariant under the location shift and
orthogonal rotation. This means that
we can assume without loss of generality that $\mathrm{E}(X_{ij})=0$ in
the rest of the paper. As noted by a reviewer, the computation of
$T_{n_1, n_2}$ would be extremely heavy if the sample sizes $n_h$ are
very large.
Indeed, the computation burden comes from the last two sums in
$A_{n_h}$ and the last three in $C_{n_1, n_2}$, where the numbers of
terms in the summations are in the order of $n_h^3$ or $n_h^4$,
respectively. Although the main motivation was the ``large $p$ small $n$'' situations, we nevertheless require $n_h \to\infty$ in our asymptotic
justifications. A solution to alleviate the computation burden can be
found by noting that the last two terms in $A_{n h}$ and the last three
in $C_{n_1, n_2}$ are all of smaller order than the first, under the
assumption of $\mu_h =0$. This means that we can first transform each
datum $X_{h i}$ to $X_{h i} - \bar{X}_{n_h}$, and then compute only
the first term in (2.1) and (2.2). These will reduce the computation to
$O(n_h^2)$ without affecting the asymptotic normality. The only price
paid for such an operation is that the modified statistic is no longer unbiased.

To establish the limiting distribution of $T_{n_1, n_2}$ so as to
establish the two sample test for the variance--covariance, we assume
the following conditions:

\begin{enumerate}[A3.]
\item[A1.] As $\min \{n_1,n_2\}\to\infty$,
$n_1/(n_1+n_2) \to\rho$ for a fixed constant $\rho\in(0,1)$.

\item[A2.] As $\min \{n_1,n_2\}\to\infty$,
$p=p(n_1,n_2) \to\infty$, and for any $k$ and $l$ $\in\{1, 2\}$,
$\operatorname{tr}(\Sigma_k\Sigma_l)\to\infty$ and
%
\begin{eqnarray}\label{eq:keyA2}
\operatorname{tr}\{(\Sigma_i\Sigma_j)(\Sigma_k\Sigma_l)\}=o\{\operatorname{tr}(\Sigma_i\Sigma
_j)\operatorname{tr}(\Sigma_k\Sigma_l)\}.
\end{eqnarray}

\item[A3.] For each $i=1$ or $2$, $X_{ij}= \Gamma_i Z_{ij}+\mu_i$
where $\Gamma_i$ is a $p \times m_i$
matrix such that $\Gamma_i \Gamma_i^{\prime}=\Sigma_i$, $\{Z_{ij}\}
_{j=1}^{n_i}$ are independent and identically
distributed (i.i.d.) $m_i$-dimensional random vectors with $m_i \ge p$
and satisfy $\mathrm{E}(Z_{ij})=0$,
$\operatorname{Var}(Z_{ij})=I_{m_i}$, the $m_i \times m_i$ identity matrix.
Furthermore, if write $Z_{ij}=(z_{ij1},\ldots ,z_{ijm_i})^{\prime}$, then
each $z_{ijk}$ has finite 8th moment,
$\mathrm{E}(z^4_{ijk})=3+\Delta_i$ for some constant $\Delta_i$ and
for any positive integers $q$ and $\alpha_l$ such that $\sum_{l=1}^q
\alpha_l \le8 \mathrm{E}(z^{\alpha_1}_{ijl_1}\cdots z^{\alpha
_q}_{ijl_q})=\mathrm{E}(z^{\alpha_1}_{ijl_1})
\cdots \mathrm{E}(z^{\alpha_q}_{ijl_q})$ for any $l_1 \ne l_2 \ne\cdots \ne l_q$.\vadjust{\goodbreak}
\end{enumerate}

While Condition A1 is of standard for two-sample asymptotic analysis,
A2~spells the extent of high dimensionality and the dependence which
can be accommodated by the proposed tests. A key aspect is that it does
not impose any explicit relationships between $p$ and the sample sizes,
but rather requires a quite mild (\ref{eq:keyA2}) regarding the covariances.
To appreciate (\ref{eq:keyA2}), we note that if $i=j=k=l$, it has the
form of $\operatorname{tr}(\Sigma_i^4)=o\{\operatorname{tr}^2(\Sigma_i^2)\}$, which is
valid if all the eigenvalues of $\Sigma_i$ are uniformly bounded.
Condition (\ref{eq:keyA2}) also makes the asymptotic study of the test
statistic manageable under high dimensionality. We note here that
requiring $\operatorname{tr}(\Sigma_k \Sigma_l) \to\infty$ is a precursor to (\ref{eq:keyA2}).
We do not assume specific parametric distributions for the two samples.
Instead, a general multivariate model
is assumed in A3 which was advocated in \citet{r5} for
testing high dimensional means.
The model resembles that of the factor model with $Z_i$ representing
the factors, except that here we allow the number of factor $m_i$ at
least as large as~$p$. This provides flexibility in accommodating a
wider range of multivariate distributions for the observed data $X_{ij}$.

Derivations leading to (\ref{stat1-final}) in Section \ref{sec6} show that,
under A2 and A3, the leading order variance of $T_{n_1,n_2}$ under
either $H_{0 a}$ or $H_{1 a}$ is
%
\begin{eqnarray}\label{variance}
\sigma^2_{n_1,n_2}
&=&\sum_{i=1}^2 \biggl[\frac{4}{n_i^2}\operatorname{tr}^2(\Sigma_i^2)
+\frac{8}{n_i}\operatorname{tr}\{(\Sigma_i^2-\Sigma_1\Sigma_2)^2\}\nonumber\\
&&\hspace*{17pt}{}+\frac{4\Delta_i}{n_i}\operatorname{tr}\{\Gamma_i^{\prime}(\Sigma_1-\Sigma
_2)\Gamma_i
\circ\Gamma_i^{\prime}(\Sigma_1-\Sigma_2)\Gamma_i\}\biggr]\\
&&{}+\frac{8}{n_1n_2}\operatorname{tr}^2(\Sigma_1 \Sigma_2),\nonumber
\end{eqnarray}
where $A \circ B = (a_{ij}b_{ij})$ for two matrices $A=(a_{ij})$ and
$B=(b_{ij})$. Note that for any symmetric matrix $A$, 
$\operatorname{tr}(A\circ A) 
\le
\operatorname{tr}(A^2)$. Hence,
\begin{eqnarray*}
\operatorname{tr}\{\Gamma_1^{\prime}(\Sigma_1-\Sigma_2)\Gamma_1 \circ\Gamma
_1^{\prime}(\Sigma_1-\Sigma_2)\Gamma_1 \}&  \le & \operatorname{tr}\{(\Sigma
_1^2-\Sigma_1\Sigma_2)^2 \}\quad  \mbox{and}\\
\operatorname{tr}\{\Gamma_2^{\prime}(\Sigma_1-\Sigma_2)\Gamma_2
\circ\Gamma_2^{\prime}(\Sigma_1-\Sigma_2)\Gamma_2\}&  \le& \operatorname{tr}\{
(\Sigma_2^2-\Sigma_2\Sigma_1)^2\}.
\end{eqnarray*}
These together with the fact that $\Delta_i \ge-2$ ensure that
$\sigma^2_{n_1,n_2} > 0$. We note that the $\Gamma_i$--$Z_{ij}$ pair
in Model A3 is not unique, and
there are other pairs, say $\tilde{\Gamma}_i$ and $\tilde{Z}_{ij}$,
such that
$X_{ij}=\tilde{\Gamma}_i \tilde{Z}_{ij}$. However, it can be shown
that the value of
$\frac{4\Delta_i}{n_i}\operatorname{tr}\{\Gamma_i^{\prime}(\Sigma_1-\Sigma
_2)\Gamma_i \circ\Gamma_i^{\prime}(\Sigma_1-\Sigma_2)\Gamma_i\}$
remains the same.

The following theorem establishes the asymptotic normality of $T_{n_1,n_2}$.

\begin{theorem}\label{teo1}
Under Conditions \textup{A1--A3}, as $\min\{n_1, n_2\}
\to\infty$
\[
\sigma_{n_1,n_2}^{-1}  [ T_{n_1,n_2}-\operatorname{tr}\{(\Sigma_1-\Sigma
_2)^2\}  ]
\stackrel{d}{\rightarrow} \mathrm{N}(0,1).
\]
\end{theorem}

It is noted that under $H_{0 a}\dvtx \Sigma_1=\Sigma_2=\Sigma$, say,
$\sigma^2_{n_1,n_2}$ becomes
\[
\sigma^2_{0,n_1,n_2}=4\biggl(\frac{1}{n_1}+\frac{1}{n_2}\biggr)^2\operatorname{tr}^2(\Sigma^2).
\]
To formulate a test procedure, we need to estimate $\sigma^2_{0, n_1,n_2}$.
As $A_{n_1}$ and $A_{n_2}$ are unbiased estimators of $\operatorname{tr}(\Sigma_1^2)$
and $\operatorname{tr}(\Sigma_2^2)$, respectively, we will use $\hat{\sigma
}^2_{0,n_1,n_2}=:\frac{2}{n_2}A_{n_1}+\frac{2}{n_1}A_{n_2}$ as the
estimator. The following theorem shows that $\hat{\sigma
}^2_{0,n_1,n_2}$ is ratio-consistent to $\sigma^2_{0,n_1,n_2}$.

\begin{theorem}\label{teo2}
Under Conditions \textup{A1--A3} and $H_{0 a}$,
as $\min\{n_1, n_2\} \to\infty$,
%
\begin{eqnarray}
\frac{A_{n_i}}{\operatorname{tr}(\Sigma_i^2)}\stackrel{p}{\rightarrow} 1\qquad  \mbox{for }
i=1\mbox{ and }2\quad  \mbox{and}\quad
\frac{\hat{\sigma}_{0,n_1,n_2}}{\sigma_{0,n_1,n_2}} \stackrel{p}{\rightarrow} 1.
\end{eqnarray}
\end{theorem}

Applying Theorems \ref{teo1} and \ref{teo2}, under $H_{0 a}\dvtx \Sigma_1=\Sigma_2$,
%
\begin{eqnarray}\label{null-stat}
L_n= \frac{T_{n_1,n_2}}{\hat{\sigma}_{0, n_1,n_2}}\stackrel{d}{\rightarrow}
\mathrm{N}(0,1).
\end{eqnarray}
Hence, the proposed test with a nominal $\alpha$ level of significance
rejects $H_{0 a}$ if
$T_{n_1,n_2} \ge\hat{\sigma}_{0, n_1,n_2} z_{\alpha}$, where
$z_{\alpha} $ is the upper-$\alpha$ quantile of N($0, 1$).

Let $\beta_{1,n_1,n_2}(\Sigma_1,\Sigma_2 ; \alpha
)=P(T_{n_1,n_2}/\hat{\sigma}_{0, n_1,n_2} >z_{\alpha}|H_{1a})$ be the
power of the test under $H_{1a}\dvtx \Sigma_1 \ne\Sigma_2$.
From Theorems \ref{teo1} and \ref{teo2}, the leading order power 
is
%
\begin{eqnarray}\label{eq:power}
\Phi\biggl( -\mathscr{Z}_{n_1,n_2}(\Sigma_1,\Sigma_2)z_{\alpha}
+\frac{\operatorname{tr}\{(\Sigma_1-\Sigma_2)^2
\}}{\sigma_{n_1,n_2}}\biggr),
\end{eqnarray}
where $\mathscr{Z}_{n_1,n_2}(\Sigma_1,\Sigma_2)= (\sigma_{n_1,n_2})^{-1}
\{\frac{2}{n_2}\operatorname{tr}(\Sigma_1^2)+\frac{2}{n_1}\operatorname{tr}(\Sigma_2^2)\}$.
It is the case that $\mathscr{Z}_{n_1,n_2}(\Sigma_1,\Sigma_2)$ is
bounded. To appreciate this, we note that $\sigma_{n_1,n_2}^2 \ge
\frac{4}{n_1^2}\operatorname{tr}^2(\Sigma_1^2)+\frac{4}{n_2^2}\operatorname{tr}^2(\Sigma_2^2)$.
Let $\gamma_p=\operatorname{tr}(\Sigma_1^2)/\operatorname{tr}(\Sigma_2^2)$ and
$k_n=n_1/(n_1+n_2)$, then
\[
\mathscr{Z}_{n_1,n_2}(\Sigma_1,\Sigma_2) \le
\frac{({2}/{n_2})\operatorname{tr}(\Sigma_1^2)+({2}/{n_1})\operatorname{tr}(\Sigma
_2^2)}{\sqrt{({4}/{n_1^2})
\operatorname{tr}^2(\Sigma_1^2)+({4}/{n_2^2})\operatorname{tr}^2(\Sigma_2^2)}} =: R_n(\gamma
_p),
\]
where $R_n(u) = (\frac{k_n}{1-k_n} u +1)\lbrace u^2+(\frac
{k_n}{1-k_n})^2\rbrace^{-1/2}$.
Since $R_n(u)$ is maximized uniquely at $u^{\ast}=(\frac
{k_n}{1-k_n})^3$, $\mathscr{Z}_{n_1,n_2}(\Sigma_1,\Sigma_2)\le\frac
{1}{k_n(1-k_n)}$. Thus,
%
\begin{eqnarray}\label{eq:lowerb}
\beta_{1,n_1,n_2}(\Sigma_1,\Sigma_2 ; \alpha)\ge
\Phi\biggl(-\frac{z_{\alpha}}{k_n(1-k_n)}
+\frac{\operatorname{tr}\{(\Sigma_1-\Sigma_2)^2
\}}{\sigma_{n_1,n_2}}\biggr)
\end{eqnarray}
implying the power is bounded from below by the probability on the
right-hand side.

Both (\ref{eq:power}) and (\ref{eq:lowerb}) indicate that $\operatorname{SNR}_1(\Sigma_1,
\Sigma_2)=:\operatorname{tr}\{(\Sigma_1-\Sigma_2)^2
\}/\sigma_{n_1,n_2}$ is instrumental in determining the power of the test.
We term $\operatorname{SNR}_1(\Sigma_1, \Sigma_2)$ as the signal-to-noise
ratio for the current testing problem since
$\operatorname{tr}\{(\Sigma_1-\Sigma_2)^2\}$ may\vadjust{\goodbreak} be viewed as the signal while
$\sigma_{n_1,n_2}$ may be viewed as the level
of the noise. If the signal is strong or the noise is weak so that the
signal-to-noise ratio diverges to the
infinity, the power will converge to 1. If the signal-to-noise ratio
diminishes to 0, the test will not be powerful
and cannot distinguish $H_{0 a}$ from $H_{1 a}$. We note that
\begin{eqnarray*}
\sigma^2_{n_1,n_2} &\le& 4\biggl\{\frac{1}{n_1}\operatorname{tr}(\Sigma_1^2)+\frac
{1}{n_2}\operatorname{tr}(\Sigma_2^2)\biggr\}^2 \nonumber\\
&&{} + \max\{8+4\Delta
_1, 8+4\Delta_2\}\biggl\{\frac{1}{n_1}\operatorname{tr}(\Sigma_1^2)+\frac
{1}{n_2}\operatorname{tr}(\Sigma_2^2)\biggr\} \operatorname{tr}\{(\Sigma_1-\Sigma_2)^2\}.
\end{eqnarray*}

Let $\delta_{1,n}=\{\frac{1}{n_1}\operatorname{tr}(\Sigma_1^2)+\frac
{1}{n_2}\operatorname{tr}(\Sigma_2^2)\}/\operatorname{tr}\{(\Sigma_1-\Sigma_2)^2\}$, then
\[
\operatorname{SNR}_1(\Sigma_1, \Sigma_2) \ge[ 4\delta_{1,n}^2+\max\{ 8+4\Delta_1, 8+4\Delta_2\}\delta_{1,n}]^{-{1}/{2}}
.
\]
Thus, if the difference between $\Sigma_1$ and $\Sigma_2$ is not too
small so that
%
\begin{equation}\label{eq:power1}
\begin{tabular}{@{}p{220pt}@{}}
$\operatorname{tr}\{(\Sigma_1-\Sigma_2)^2\}$ is at the same or a larger
order of $\frac{1}{n_1}\operatorname{tr}(\Sigma
_1^2)+\frac{1}{n_2}\operatorname{tr}(\Sigma_2^2)$,
\end{tabular}
\end{equation}
the test will be powerful. Condition (\ref{eq:power1}) is trivially
true for fixed-dimen\-sional data while $n_i \to\infty$.
For high-dimensional data, it is less automatic as $\operatorname{tr}(\Sigma_i^2)$
can diverge.
To gain further insight on (\ref{eq:power1}), let $\lambda_{i1} \le
\lambda_{i2} \le\cdots\le\lambda_{ip}$ be
the eigenvalues of $\Sigma_i$. Then, a sufficient condition for the
test to
have a~nontrivial power is $\operatorname{tr}\{(\Sigma_1-\Sigma_2)^2\}= O\{\frac
{1}{n_1}\sum_{i=1}^{p}\lambda_{1i}^2+\frac{1}{n_2}
\sum_{i=1}^{p}\lambda_{2i}^2\}$. If all the eigenvalues of $\Sigma
_1$ and $\Sigma_2$ are bounded away from zero and infinity,~(\ref{eq:power1}) becomes $\operatorname{tr}\{(\Sigma_1-\Sigma_2)^2\}=O(n^{-1}p)$.
Let $\delta_{\beta}= p^{-1} \sqrt{\operatorname{tr}\{(\Sigma_1-\Sigma_2)^2\}}$ be
the average signal. Then the test has nontrivial power if $\delta
_{\beta}$ is at least at the order of $n^{-{1}/{2}}p^{-{1}/{2}}$, which is actually smaller than the conventional order of
$n^{-1/2}$ for fixed-dimension situations. This partially reflects the
fact that high data dimensionality is not entirely a curse as
there are more data information available as well. 
If the covariance matrix is believed to have certain structure, for
instance banded or bandable in the sense of \citet{r9},
we may modify the test statistic so that the comparison of the two
covariance matrices is made in the ``important regions'' under the
structure. The modification can be in the form of thresholding, a topic
we would not elaborate in this paper; see \citet{r12} for
research in this direction.

\section{Test for covariance between two sub-vectors}\label{sec3}

Let $X_{i j}=(X_{ij}^{(1)}, X_{ij}^{(2)})$ be\vspace*{2pt} a partition of the
original data vector into sub-vectors of dimensions of $p_1$ and~$p_2$,
and $\Sigma_{i, 12} =\operatorname{Cov}(X_{ij}^{(1)}, X_{ij}^{(2)})$ be the
covariance between the sub-vectors. The\vspace*{1pt} focus in this section is to
develop a test procedure for $H_{0 b}\dvtx \Sigma_{1, 12} = \Sigma_{2,
12}$. Testing for such a hypothesis is importance in its own right, for
instance in detecting changes in correlation between two groups of\vadjust{\goodbreak}
genes under two treatment regimes. It can be also viewed as part of the
effort in reducing the dimensionality in testing high-dimensional
variance--covariances. To elaborate on this, consider the partition of~$\Sigma_i$,
%
\begin{equation}\label{eq:block}
\Sigma_i = \pmatrix{
 \Sigma_{i, 11} & \Sigma_{i, 12} \cr
\Sigma_{i, 12}^{\prime} & \Sigma_{i, 22}
},
\end{equation}
induced by the partition of the data vectors. Instead of testing on the
whole matrices
$\Sigma_1 = \Sigma_2$, we can first test separately on the two
diagonal blocks $\Sigma_{1, ll} = \Sigma_{2, ll}$ for $l=1$ and $2$,
by employing the test developed in the previous section based on the
sub-vectors of the two sample data respectively. Then, we can test for
the off-diagonal blocks $H_{0 b}\dvtx \Sigma_{1, 12} = \Sigma_{2, 12}$
using a test procedure to be developed in this section. 


The partition of data vectors also induces a partition of the
multivariate model in A3 so that
%
\begin{eqnarray}\label{stat2-6}
X_{ij}^{(1)}= \Gamma_i^{(1)} Z_{ij}+\mu_i^{(1)}\quad  \mbox{and}
\quad
X_{ij}^{(2)}= \Gamma_i^{(2)} Z_{ij}+\mu_i^{(2)},
\end{eqnarray}
where $\Gamma_i^{(1)}$ is $p_1 \times m_i$ and $\Gamma_i^{(2)}$ is
$p_2 \times m_i$ such that $\Gamma_i^{\prime}=({\Gamma
_i^{(1)\prime}}, {\Gamma_i^{(2)\prime}})$ and $\Gamma
_i^{(1)}{\Gamma_i^{(2)\prime}}=\Sigma_{i,12}$.

We are interested in testing $H_{0 b}\dvtx \Sigma_{1,12}=\Sigma_{2,12}$
vs $H_{1 b}\dvtx \Sigma_{1,12} \ne\Sigma_{2,12}$. The test statistic is
aimed at
\begin{eqnarray}\label{stat2-1}
& &\operatorname{tr}\{(\Sigma_{1,12}-\Sigma_{2,12})(\Sigma_{1,12}-\Sigma
_{2,12})^{\prime} \}\nonumber\\[-9.5pt]\\[-9.5pt]
&&\qquad=\operatorname{tr}(\Sigma_{1,12}\Sigma_{1,12}^{\prime})+\operatorname{tr}(\Sigma_{2,12}\Sigma
_{2,12}^{\prime})
-2\operatorname{tr}(\Sigma_{1,12}\Sigma_{2,12}^{\prime}),\nonumber
\end{eqnarray}
a discrepancy measure between $\Sigma_{1,12}$ and $\Sigma_{2,12}$.

With the same considerations as those when we proposed the estimators
in (\ref{eqn2}) and (\ref{eqn3}),
we estimate $\operatorname{tr}(\Sigma_{h,12}\Sigma_{h,12}^{\prime})$ by
%
\begin{eqnarray}\label{stat2-2}
U_{n_h}
&=& \frac{1}{n_h(n_h-1)} \sum_{i \ne j} {X_{hi}^{(1)\prime}}X_{hj}^{(1)}
{X_{hj}^{(2)\prime}}X_{hi}^{(2)}\nonumber\\[-4pt]
&&{}- \frac{2}{n_h(n_h-1)(n_h-2)} \sum_{i,j,k}^{\star}
{X_{hi}^{(1)\prime}} X_{hj}^{(1)}{X_{hj}^{(2)\prime}}
X_{hk}^{(2)}\\[-4pt]
&&{} + \frac{1}{n_h(n_h-1)(n_h-2)(n_h-3)} \sum_{i,j,k,l}^{\star}
{X_{hi}^{(1)\prime}} X_{hj}^{(1)}{X_{hk}^{(2)\prime}}
X_{hl}^{(2)}, \nonumber
\end{eqnarray}
and estimate $\operatorname{tr}(\Sigma_{1,12}\Sigma_{2,12}^{\prime})$ by
\begin{eqnarray}\label{stat2-4}
W_{n_1n_2}
&=&
\frac{1}{n_1n_2}\sum_{i,j}{X_{1i}^{(1)\prime}
}X_{2j}^{(1)}{X_{2j}^{(2)\prime}}X_{1i}^{(2)}\nonumber\\[-4pt]
&&{}- \frac{1}{n_1
n_2 (n_1-1)}\sum_{i \ne k, j}{X_{1 i}^{(1)\prime}}X_{2j}^{(1)}{X_{2
j}^{(2)\prime}}X_{1 k}^{(2)}\nonumber\\[-9.5pt]\\[-9.5pt]
&&{}- \frac{1}{n_1 n_2 (n_2-1)}\sum_{i \ne k, j}{X_{2
i}^{(1)\prime}}X_{1j}^{(1)}{X_{1 j}^{(2)\prime}}X_{2 k}^{(2)}
\nonumber \\[-4pt]
&&{}+ \frac{1}{n_1 n_2 (n_1-1)(n_2-1)}\sum_{i \ne k, j\ne l}{X_{1
i}^{(1)\prime}}X_{2j}^{(1)}{X_{1 k}^{(2)\prime}}X_{2
l}^{(2)}.\nonumber
\end{eqnarray}
Both $U_{n h}$ and $W_{n_1 n_2}$ are linear combinations of U-statistics.

Combining these estimators together leads to an unbiased estimator
of
$\operatorname{tr}\{(\Sigma_{1,12}-\Sigma_{2,12})(\Sigma_{1,12}-\Sigma
_{2,12})^{\prime} \}$,
%
\begin{eqnarray}\label{stat2-5}
S_{n_1,n_2}=U_{n_1}+U_{n_2}-2 W_{n_1n_2},
\end{eqnarray}
which is also invariant under the location shift and orthogonal rotations.

To establish the asymptotic normality
of $S_{n_1,n_2}$, we need an extra assumption regarding the
off-diagonal sub-matrices.

\begin{enumerate}[A4.]
\item[A4.] As $\min \{n_1,n_2\}\to\infty$, for any $i,j,k$
and $l$ $\in\{1, 2\}$.
%
\begin{eqnarray}
\operatorname{tr}(\Sigma_{i,11}\Sigma_{j,12}\Sigma_{k,22}\Sigma_{l,12}^{\prime})
=o\{\operatorname{tr}(\Sigma_{i,11}\Sigma_{j,11})\operatorname{tr}(\Sigma_{k,22}\Sigma_{l,22})\}.
\end{eqnarray}
\end{enumerate}

Derivations leading to (\ref{stat2-final}) in Section \ref{sec6} show that,
under A2, A3 and A4, the leading order variance of $S_{n_1,n_2}$ is
%
\begin{eqnarray}\label{stat2-variance}
\omega^2_{n_1,n_2}
&=&\sum_{i=1}^2 \biggl[\frac{2}{n_i^2}\operatorname{tr}^2(\Sigma_{i,12}\Sigma
_{i,12}^{\prime})
+\frac{2}{n_i^2}\operatorname{tr}(\Sigma_{i,11}^2)\operatorname{tr}(\Sigma_{i,22}^2)\nonumber\\
&&\hspace*{18pt}{}+\frac{4}{n_i}\operatorname{tr}\{(\Sigma_{i,12}\Sigma_{1,12}^{\prime}-\Sigma
_{i,12}\Sigma_{2,12}^{\prime})^2\}\nonumber\\
&&\hspace*{18pt}{} + \frac{4}{n_i}\operatorname{tr}\{(\Sigma_{i,11}\Sigma_{1,12}-\Sigma
_{i,11}\Sigma_{2,12})
(\Sigma_{i,22}\Sigma_{1,12}^{\prime}-\Sigma_{i,22}\Sigma
_{2,12}^{\prime})\}\\
&&\hspace*{18pt}{} + \frac{4\Delta_i}{n_i}\operatorname{tr}\bigl\{{\Gamma_i^{(1)\prime}}(\Sigma
_{1,12}-\Sigma_{2,12})\Gamma_i^{(2)}
\circ{\Gamma_i^{(1)\prime}}(\Sigma_{1,12}-\Sigma_{2,12})\Gamma
_i^{(2)}\bigr\} \biggr]\nonumber\\
&&{}+\frac{4}{n_1n_2}
\operatorname{tr}^2(\Sigma_{1,12}\Sigma_{2,12}^{\prime})
+ \frac{4}{n_1n_2}
\operatorname{tr}(\Sigma_{1,11}\Sigma_{2,11})\operatorname{tr}(\Sigma_{1,22}\Sigma_{2,22}).\nonumber
\end{eqnarray}

Similarly to the analysis on $T_{n_1,n_2}$ in the previous section, the
asymptotic normality of $S_{n_1,n_2}$ can be established in the
following theorem.

\begin{theorem}\label{teo3}
 Under Conditions \textup{A1--A4}, as $\min \{
n_1,n_2\}\to\infty$,
\[
{\omega_{n_1,n_2}}^{-1}  [ S_{n_1,n_2}-\operatorname{tr}\{(\Sigma
_{1,12}-\Sigma_{2,12})(\Sigma_{1,12}-\Sigma_{2,12})^{\prime}  ]
\stackrel{d}{\rightarrow}  \mathrm{N}(0,1).
\]
\end{theorem}

Under $H_{0 b}\dvtx\Sigma_{1,12}=\Sigma_{2,12}=\Sigma_{12}$, say,
$\omega^2_{n_1,n_2}$ becomes
\begin{eqnarray}\qquad
\omega^2_{0,n_1,n_2}
&=&2\biggl(\frac{1}{n_1}+\frac{1}{n_2}\biggr)^2\operatorname{tr}^2(\Sigma_{12}\Sigma
_{12}^{\prime}) + 2 \sum_{i=1}^2 \frac{1}{n_i^{2}} \operatorname{tr}(\Sigma
_{i,11}^2)\operatorname{tr}(\Sigma_{i,22}^2) \nonumber\\[-9.5pt]\\[-9.5pt]
&&{}+ \frac{4}{n_1n_2}
\operatorname{tr}(\Sigma_{1,11}\Sigma_{2,11})\operatorname{tr}(\Sigma_{1,22}\Sigma_{2,22}).\nonumber\vadjust{\goodbreak}
\end{eqnarray}

In order to formulate a test procedure, $\omega^2_{0, n_1,n_2}$ needs
to be estimated. An unbiased estimator of $\operatorname{tr}(\Sigma_{h,ll}^2)$ for
$h=1$ or $2$ and $l=1$ or $2$, is
\begin{eqnarray*}
A_{n_h}^{(l)}
&=& \frac{1}{n_h(n_h-1)} \sum_{i \ne j} \bigl({X_{hi}^{(l)\prime}}
X_{hj}^{(l)}\bigr)^2-\frac{2}{n_h(n_h-1)(n_h-2)} \sum_{i,j,k}^{\star}
{X_{hi}^{(l)\prime}} X_{hj}^{(l)}{X_{hj}^{(l)\prime}}
X_{hk}^{(l)}\\
& &{}+ \frac{1}{n_h(n_h-1)(n_h-2)(n_h-3)} \sum_{i,j,k,l}^{\star}
{X_{hi}^{(l)\prime}} X_{hj}^{(l)}{X_{hk}^{(l)\prime}}
X_{hl}^{(l)}.
\end{eqnarray*}

Similarly, an unbiased estimator of $\operatorname{tr}(\Sigma_{1,hh}\Sigma_{2,hh})$,
for $h=1$ or $2$, is
\begin{eqnarray*}
C_{n_1n_2}^{(h)}
&=&
\frac{1}{n_1n_2}\sum_{i,j}\bigl({X_{1i}^{(h)\prime}}X_{2j}^{(h)}\bigr)^2
-\frac{1}{n_1
n_2 (n_1-1)}\sum_{i \ne k, j}{X_{1 i}^{(h)\prime}}X_{2j}^{(h)}{X_{2
j}^{(h)\prime}}X_{1 k}^{(h)}\\
&&{}- \frac{1}{n_1 n_2 (n_2-1)}\sum_{i \ne k, j}{X_{2
i}^{(h)\prime}}X_{1j}^{(h)}{X_{1 j}^{(h)\prime}}X_{2
k}^{(h)}\\
&&{} + \frac{1}{n_1 n_2 (n_1-1)(n_2-1)}\sum_{i \ne k, j\ne l}{X_{1
i}^{(h)\prime}}X_{2j}^{(h)}{X_{1 k}^{(h)\prime}}X_{2
l}^{(h)}.
\end{eqnarray*}
Then under $H_{0 b}$, an unbiased estimator of $\omega^2_{0, n_1,n_2}$ is
\[
\widehat{\omega}^2_{0, n_1,n_2}
=2\biggl(\frac{U_{n_1}}{n_2}+\frac{U_{n_2}}{n_1}\biggr)^2+\frac
{2}{n_1^2}A_{n_1}^{(1)}A_{n_1}^{(2)}+\frac{2}{n_2^2}A_{n_2}^{(1)}A_{n_2}^{(2)}
+\frac{4}{n_1n_2}C_{n_1n_2}^{(1)}C_{n_1n_2}^{(2)}.
\]

The following theorem shows that $\widehat{\omega}^2_{0, n_1,n_2}$ is
ratio-consistent to $\omega^2_{0, n_1,n_2}$.

\begin{theorem}\label{teo4}
 Under Conditions \textup{A1--A4}, and $H_{0 b}\dvtx \Sigma
_{1,12}=\Sigma_{2,12}$,
\[
\frac{\widehat{\omega}^2_{0, n_1,n_2}}{\omega^2_{0,
n_1,n_2}}\stackrel{p}{\rightarrow}  1.
\]
\end{theorem}

Applying Theorems \ref{teo3} and \ref{teo4}, we have, under $H_{0 b}$,
\[
\frac{S_{n_1,n_2}}{\hat{\omega}_{0, n_1,n_2}}\stackrel{d}{\rightarrow}
\mathrm{N}(0,1).
\]
This suggests an $\alpha$-level test that rejects $H_{0 b}$ if
$S_{n_1,n_2} \ge\hat{\omega}_{0, n_1,n_2} z_{\alpha}$. The power of
the proposed test under $H_{1 b}\dvtx \Sigma_{1,12} \ne\Sigma_{2,12}$ is
\[
\beta_{2,n_1,n_2}(\Sigma_{1,12},\Sigma_{2,12} ; \alpha)
= P(S_{n_1,n_2}/\hat{\omega}_{0, n_1,n_2} >z_{\alpha
}|H_{1b}).
\]
From Theorems \ref{teo3} and \ref{teo4}, the leading order power is
\[
\Phi\biggl(-\frac{\tilde{\omega}}
{\omega_{n_1,n_2}}z_{\alpha}
+\frac{\operatorname{tr}\{(\Sigma_{1,12}-\Sigma_{2,12})(\Sigma_{1,12}-\Sigma
_{2,12})^{\prime}
\}}{\omega_{n_1,n_2}}\biggr),
\]
where
\begin{eqnarray*}
\tilde{\omega}^2& = & 2\biggl\{\frac{\operatorname{tr}(\Sigma_{1,12}\Sigma
_{1,12}^{\prime})}{n_2} +\frac{\operatorname{tr}(\Sigma_{2,12}\Sigma
_{2,12}^{\prime})}{n_1}\biggr\}^2+ \frac{2}{n_1^2}\operatorname{tr}(\Sigma_{1,11}^2)
\operatorname{tr}(\Sigma_{1,22}^2)  \\
&&{}+  \frac{2}{n_2^2}\operatorname{tr}(\Sigma_{2,11}^2) \operatorname{tr}(\Sigma_{2,22}^2)+\frac
{4}{n_1n_2}\operatorname{tr}(\Sigma_{1,11}\Sigma_{2,11})
\operatorname{tr}(\Sigma_{1,22}\Sigma_{2,22}).
\end{eqnarray*}
Let $\eta_p=\operatorname{tr}(\Sigma_{1,12}\Sigma_{1,12}^{\prime})/\operatorname{tr}(\Sigma
_{2,12}\Sigma_{2,12}^{\prime})$.
It may be shown that 
%
\[
\frac{\tilde{\omega}}
{\omega_{n_1,n_2}}
\le
\sqrt{R^2(\eta_p)+1},
\]
where $R(\gamma_p)$ is the same function defined in Section \ref{sec2}.
Hence, asymptotically,
\begin{eqnarray*}
& &\beta_{2,n_1,n_2}(\Sigma_{1,12},\Sigma_{2,12} ; \alpha
)\\
&&\qquad\ge
\Phi\biggl(-\frac{z_{\alpha}\sqrt{1+k_n^2(1-k_n)^2}}{k_n(1-k_n)}
+\frac{\operatorname{tr}\{(\Sigma_{1,12}-\Sigma_{2,12})(\Sigma_{1,12}-\Sigma
_{2,12})^{\prime}
\}}{\omega_{n_1,n_2}}\biggr).
\end{eqnarray*}
This implies that
\[
\operatorname{SNR}_2 =: \operatorname{tr}\{(\Sigma_{1,12}-\Sigma_{2,12})(\Sigma
_{1,12}-\Sigma_{2,12})^{\prime}
\}/\omega_{n_1,n_2}
\]
is the key quantity that determines the power of the test. Furthermore, let
%
\[
\delta_{2,n}=\frac{({1}/{n_1})\operatorname{tr}(\Sigma_{1,11})\operatorname{tr}(\Sigma
_{1,22})+({1}/{n_2})\operatorname{tr}(\Sigma_{2,11})\operatorname{tr}(\Sigma_{2,22})}
{\operatorname{tr}\{(\Sigma_{1,12}-\Sigma_{2,12})(\Sigma_{1,12}-\Sigma
_{2,12})^{\prime} \}}.
\]
%
It can be shown that
%
\begin{eqnarray}\label{eq:snr2}
\operatorname{SNR}_2\ge [4\delta_{2,n}^2+\max\{8+4\Delta_1,
8+4\Delta_2\}\delta_{2,n} ]^{-{1}/{2}}
.
\end{eqnarray}
Hence, the test is powerful if the difference between $\Sigma_{1,12}$
and $\Sigma_{2,12}$ is not too small so that
$\operatorname{tr}\{(\Sigma_{1,12}-\Sigma_{2,12})(\Sigma_{1,12}-\Sigma
_{2,12})^{\prime}\}$ is at the order of
$\sum_{i=1}^2 \frac{1}{n_i}\times\operatorname{tr}(\Sigma_{i,11})\operatorname{tr}(\Sigma_{i,22})$ or larger.
A further analysis on the power, similar to that given at the end of
last section, can be made. Here for the sake of brevity, we will not report.

\section{Simulation studies}\label{sec4}

We report results from simulation experiments which were designed to
evaluate the performance of the two proposed
tests. A range of dimensionality and sample sizes was considered which
allowed $p$ to increase as the sample sizes were increased. This was
designed to confirm the asymptotic results reported in the previous sections.

We first considered the test for $H_{0 a}\dvtx \Sigma_1 = \Sigma_2$
regarding the whole variance--covariance matrices. To compare with the
conventional likelihood ratio (LR) test and the corrected LR test
proposed by \citet{r4}, we first considered cases of $p \le\min
\{n_1, n_2\}$ and the normally distributed data.
Specifically, to create the null hypothesis, we simulated both samples
from the $p$-dimensional standard normal distribution.
To evaluate the power of the three tests, we set the first population
to be the $p$-dimensional standard normally distributed while
simulating the second population according to 
%
\begin{eqnarray}\label{eq:ma1}
X_{ijk}=Z_{ijk}+\theta_1 Z_{ijk+1},
\end{eqnarray}
where $\{Z_{ijk}\}$ were i.i.d. standard normally distributed, and
$\theta_1 =0.5, 0.3$ and $0.2$, respectively. As $\theta_1$ was
decreased, the signal strength for the test became weaker.
We chose $(p,n_1,n_2)=(40, 60, 60), (80, 120, 120)$ and $(120, 180,
180)$, respectively. The empirical size and power for the three tests
are reported in Table \ref{tab1}. All the simulation results reported in this
section were based on 1000 simulations with the nominal significance
level to be 5$\%$.

\begin{table}
\caption{Empirical sizes and powers of the conventional likelihood
ratio (LR), the corrected likelihood ratio (CLR) and the proposed tests
(Proposed) for the variance--covariance, based on 1000 replications with
normally distributed $\{Z_{ijk}\}$} \label{tab1}
\begin{tabular*}{\textwidth}{@{\extracolsep{\fill}}lcd{1.3}d{1.3}d{1.3}d{1.3}@{}}
\hline
&&&\multicolumn{3}{c@{}}{\textbf{Power}}\\[-5pt]
&&&\multicolumn{3}{c@{}}{\hrulefill}\\
$\bolds{(p,n_1,n_2)}$& \textbf{Methods}& \multicolumn{1}{c}{\textbf{Size}}& \multicolumn{1}{c}{$\bolds{\theta_1=0.5}$}
&\multicolumn{1}{c}{$\bolds{\theta_1=0.3}$}& \multicolumn{1}{c@{}}{$\bolds{\theta_1=0.2}$}\\
\hline
$(40,60,60)$ & LRT & 1 &  1 &  1 & 1 \\
& CLRT& 0.043&  0.999 &  0.509& 0.172\\
& Proposed & 0.052&  0.999 &  0.734& 0.271\\
$(80,120,120)$& LRT & 1 & 1 &  1 &  1 \\
& CLRT& 0.045& 1 &  0.946&  0.421\\
& Proposed & 0.053& 1 &  0.997& 0.713\\
$(120,180,180)$& LRT& 1 & 1 &  1 & 1 \\
& CLRT& 0.062& 1 &  1 & 0.713 \\
& Proposed & 0.045& 1 &  1 & 0.958\\
\hline
\end{tabular*}
\end{table}

We then carried out simulations for situations where $p$ was much
larger than the sample sizes. In this case,
only the proposed test was considered, as both the LR and the corrected
LR tests were no longer applicable.
We chose a set of data dimensions from $32$ to $700$, while the sample
sizes ranged from 20 to 100, respectively.
We considered the moving average model (\ref{eq:ma1}) with $\theta
_1=2$ as the null model of both populations for size evaluation.
To assess the power performance, the first population was generated
according to (\ref{eq:ma1}), while the second was from
%
\begin{eqnarray}\label{eq:ma2}
X_{ijk}=Z_{ijk}+\theta_1 Z_{ijk+1}+\theta_2 Z_{ijk+2},
\end{eqnarray}
where $\theta_1=2$ and $\theta_2=1$.
Three combinations of distributions were experimented for the i.i.d.
sequences $\{Z_{ijk}\}_{k=1}^p$ in models (\ref{eq:ma1}) and (\ref{eq:ma2}), respectively. They were: (i) both sequences were the
standard normal; (ii) the centralized $\operatorname{Gamma}(4,0.5)$ for Sample 1 and the
centralized $\operatorname{Gamma}(0.5,\sqrt{2}$) for Sample 2; (iii) the standard
normal for Sample 1 and the centralized $\operatorname{Gamma}(0.5,\sqrt{2}$) for
Sample 2. The last two combinations were designed to assess the
performance under nonnormality. The empirical size and power of the
test are reported in Tables~\ref{tab2}--\ref{tab4}.

\begin{table}
\caption{Empirical sizes and powers of the proposed test for the
variance--covariance matrices, based on 1000 replications with normally
distributed $\{Z_{ijk}\}$ in Models (\protect\ref{eq:ma1}) and (\protect\ref
{eq:ma2})} \label{tab2}
\begin{tabular*}{\textwidth}{@{\extracolsep{\fill}}lcccccc@{}}
\hline
&\multicolumn{6}{c@{}}{$\bolds{p}$}\\[-5pt]
&\multicolumn{6}{c@{}}{\hrulefill}\\
$\bolds{n_1=n_2}$&\multicolumn{1}{c}{$\bolds{32}$} & \multicolumn{1}{c}{$\bolds{64}$} &
\multicolumn{1}{c}{$\bolds{128}$} &\multicolumn{1}{c}{$\bolds{256}$}
&\multicolumn{1}{c}{$\bolds{512}$}&\multicolumn{1}{c@{}}{$\bolds{700}$}\\
\hline
\multicolumn{7}{@{}c@{}}{Sizes}\\
\phantom{0}20 & 0.044 & 0.054 &0.051 & 0.048 & 0.051 & 0.038\\
\phantom{0}50 & 0.052 & 0.060 &0.033 & 0.043 & 0.054 & 0.049\\
\phantom{0}80 & 0.054 & 0.060 &0.047 & 0.048 & 0.052 & 0.053\\
100 & 0.056 & 0.049 &0.052 & 0.046 & 0.049 & 0.048\\[6pt]
\multicolumn{7}{@{}c@{}}{Powers}\\
\phantom{0}20 & 0.291 & 0.256 & 0.267 & 0.277 & 0.282 & 0.291\\
\phantom{0}50 & 0.746 & 0.821 & 0.830 & 0.837 & 0.832 & 0.849\\
\phantom{0}80 & 0.957 & 0.992 & 0.991 & 0.998 & 0.999 & 0.998\\
100& 0.994 & 1\phantom{000.} & 0.999 & 1\phantom{000.} & 1\phantom{000.} & 1\phantom{000.}\\
\hline
\end{tabular*}
\end{table}

\begin{table}[b]
\caption{Empirical sizes and powers of the proposed test for the
variance--covariance matrices, based on 1000 replications with Gamma
distributed $\{Z_{ijk}\}$ in Models (\protect\ref{eq:ma1}) and (\protect\ref{eq:ma2})}
\label{tab3}
\begin{tabular*}{\textwidth}{@{\extracolsep{\fill}}lcccccc@{}}
\hline
&\multicolumn{6}{c@{}}{$\bolds{p}$}\\[-5pt]
&\multicolumn{6}{c@{}}{\hrulefill}\\
$\bolds{n_1=n_2}$&\multicolumn{1}{c}{$\bolds{32}$} & \multicolumn{1}{c}{$\bolds{64}$} &
\multicolumn{1}{c}{$\bolds{128}$} &\multicolumn{1}{c}{$\bolds{256}$}
&\multicolumn{1}{c}{$\bolds{512}$}&\multicolumn{1}{c@{}}{$\bolds{700}$}\\
\hline
\multicolumn{7}{@{}c@{}}{Sizes}\\
\phantom{0}20 & 0.119 & 0.117 & 0.069 & 0.063 & 0.051 &0.040\\
\phantom{0}50 & 0.150 & 0.110 & 0.094 & 0.052 & 0.053 &0.051\\
\phantom{0}80 & 0.155 & 0.111 & 0.093 & 0.067 & 0.064 &0.044\\
100 & 0.148 & 0.120 & 0.084 & 0.056 & 0.058 &0.053\\[6pt]
\multicolumn{7}{@{}c@{}}{Powers}\\
\phantom{0}20 & 0.299 & 0.282 &0.290 & 0.309 &0.265 & 0.277\\
\phantom{0}50 & 0.574 & 0.665 &0.693 & 0.750 &0.801 & 0.828\\
\phantom{0}80 & 0.804 & 0.886 &0.942 & 0.968 &0.991 & 0.986\\
$100$& 0.899 & 0.945 &0.986 & 0.995 &0.998 & 1\phantom{000.}\\
\hline
\end{tabular*}
\end{table}

We observed from Table \ref{tab1} that the size of the conventional LR test was
grossly distorted, confirming its breakdown under even mild
dimensionality, discovered in \citet{r4}. The severely distorted
size for the LR test made its power artificially high.
Both the corrected LR test and the proposed test had quite accurate
size approximation to the nominal 5\% level for all cases in Table \ref{tab1}.
Both tests enjoyed perfect power at $\theta_1=0.5$, when the signal
strength of the tests was strong.
When the value of $\theta_2$ decreased, both tests had smaller power,
although the proposed test was slightly more powerful than the
corrected LR test at $\theta_1=0.3$ and much more so at $\theta
_1=0.2$, when the signal strength was weaker.

The simulation results for the proposed test with dimensions much
larger than the sample sizes and for nonnormally distributed data are
reported in Tables~\ref{tab2}--\ref{tab4}. We note that the LR tests are not applicable
for the setting. The simulation results show that the proposed test had
quite accurate and robust size approximation in a quite wider range of
dimensionality and distributions, considered in the simulation
experiments. The tables also show that
the power of the proposed tests was quite satisfactory and was
increased as the dimension and the sample sizes became larger.

\begin{table}
\caption{Empirical sizes and powers of the proposed test for the
variance--covariance matrices, based on 1000 replications with the mixed
normal and Gamma distributions for $\{Z_{ijk}\}$ in Models (\protect\ref
{eq:ma1}) and~(\protect\ref{eq:ma2})} \label{tab4}
\begin{tabular*}{\textwidth}{@{\extracolsep{\fill}}lcccccc@{}}
\hline
&\multicolumn{6}{c@{}}{$\bolds{p}$}\\[-5pt]
&\multicolumn{6}{c@{}}{\hrulefill}\\
$\bolds{n_1=n_2}$&\multicolumn{1}{c}{$\bolds{32}$} & \multicolumn{1}{c}{$\bolds{64}$} &
\multicolumn{1}{c}{$\bolds{128}$} &\multicolumn{1}{c}{$\bolds{256}$}
&\multicolumn{1}{c}{$\bolds{512}$}&\multicolumn{1}{c@{}}{$\bolds{700}$}\\
\hline
\multicolumn{7}{@{}c@{}}{Sizes}\\
\phantom{0}20 & 0.108 & 0.099 &0.076 & 0.059 & 0.070 & 0.050\\
\phantom{0}50 & 0.117 & 0.111 &0.069 & 0.068 & 0.057 & 0.053\\
\phantom{0}80 & 0.124 & 0.099 &0.091 & 0.065 & 0.064 & 0.060\\
100 & 0.150 & 0.122 &0.085 & 0.069 & 0.056 & 0.047\\[6pt]
\multicolumn{7}{@{}c@{}}{Powers}\\
\phantom{0}20 & 0.256 & 0.296 & 0.278 & 0.297 & 0.276 & 0.295\\
\phantom{0}50 & 0.606 & 0.659 & 0.724 & 0.766 & 0.824 & 0.823\\
\phantom{0}80 & 0.805 & 0.890 & 0.950 & 0.977 & 0.989 & 0.992\\
100 & 0.904 & 0.958 & 0.982 & 0.996 & 0.999 & 1\phantom{000.}\\
\hline
\end{tabular*}
\end{table}

We then conducted simulations to evaluate the performance of the second test
for $H_{0 b}\dvtx \Sigma_{1,12}=\Sigma_{2,12}$. We partition equally the
entire random vector~$X_{ij}$ into two subvectors of $p_1 =p/2$ and
$p_2=p-p_1$. 
To ensure sufficient number of nonzero elements in the off-diagonal
sub-matrices $\Sigma_{1,12}$ and $\Sigma_{2,12}$ when the dimension
was increased,
we considered a moving average model of order~$m_1$, which is much
larger than the orders used
in (\ref{eq:ma1}) and (\ref{eq:ma2}).
In the size evaluation,
%
\begin{eqnarray}\label{eq:stat2-ma1}
X_{ijk}=Z_{ijk}+\alpha_1 Z_{ijk+1}+\cdots+\alpha
_{m_1}Z_{ijk+m_1}
\end{eqnarray}
for $i=1,2$, $j=1,\ldots, n_i$, where all the $\alpha_i$ coefficients
were chosen to be 0.1.
In the simulation for the power, we generated the first sample
according to the above (\ref{eq:stat2-ma1}) and the second from
%
\begin{eqnarray}\label{eq:stat2-ma2}
X_{ijk}=Z_{ijk}+\beta_1 Z_{ijk+1}+\cdots+\beta_{m_2} Z_{ijk+
m_2}
\end{eqnarray}
for $j=1,\ldots, n_2$, where the $\beta_i$ were chosen to be 0.8.
We chose the lengths of the moving average $m_1$ and $m_2$ according to
the dimension $p$ such that as $p$ was increased, the values of $m_1$
and $m_2$ were increased as well. Specifically, we set $(m_1,m_2, p) =
(2, 25, 50), (3, 50, 100), (7, 100, 200), (12, 250, 500)$ and $(18,
300, 700)$, respectively.
Two distributions were considered for the i.i.d. sequences $\{Z_{ijk}\}
_{k=1}^p$ in (\ref{eq:stat2-ma1}) and (\ref{eq:stat2-ma2}): (i) both
sequences were standard normally distributed; (ii) the centralized
$\operatorname{Gamma}(4,0.5)$ for Sample 1 and the centralized $\operatorname{Gamma}(0.5,\sqrt{2}$)
for Sample 2. The simulation results for the second test are reported
in Table \ref{tab5} for the normally distributed case and Table \ref{tab6} for the Gamma
distributed case.

\begin{table}
\caption{Empirical sizes and powers of the proposed test for the
covariance between two sub-vectors, based on 1000 replications for
normally distributed $\{Z_{ijk}\}$ in Models (\protect\ref{eq:stat2-ma1}) and
(\protect\ref{eq:stat2-ma2})}
\label{tab5}
\begin{tabular*}{\textwidth}{@{\extracolsep{\fill}}lccccc@{}}
\hline
&\multicolumn{5}{c@{}}{$\bolds{p}$}\\[-5pt]
&\multicolumn{5}{c@{}}{\hrulefill}\\
 $\bolds{n_1=n_2}$&\multicolumn{1}{c}{$\bolds{50}$} &
\multicolumn{1}{c}{$\bolds{100}$} &\multicolumn{1}{c}{$\bolds{200}$}
&\multicolumn{1}{c}{$\bolds{500}$}&\multicolumn{1}{c@{}}{$\bolds{700}$} \\
\hline
\multicolumn{6}{@{}c@{}}{Sizes}\\
\phantom{0}20 & 0.069 & 0.071 & 0.070 & 0.065 & 0.077\\
\phantom{0}50 & 0.064 & 0.056 & 0.064 & 0.063 & 0.055\\
\phantom{0}80 & 0.057 & 0.046 & 0.056 & 0.073 & 0.052\\
100 & 0.047 & 0.062 & 0.055 & 0.054 & 0.048\\[6pt]
\multicolumn{6}{@{}c@{}}{Powers}\\
\phantom{0}20 & 0.639 & 0.625 & 0.628 & 0.620& 0.615\\
\phantom{0}50 & 0.993 & 0.994 & 0.982 & 0.983& 0.989\\
\phantom{0}80 & 1\phantom{000.} & 1\phantom{000.} & 1\phantom{000.} & 1\phantom{000.} & 1\phantom{000.}\\
100& 1\phantom{000.} & 1\phantom{000.} & 1\phantom{000.} & 1\phantom{000.} & 1\phantom{000.}\\
\hline
\end{tabular*}
\end{table}

\begin{table}
\caption{Empirical sizes and powers of the proposed test for the
covariances between two sub-vectors, based on 1000 replications
with Gamma distributed $\{Z_{ijk}\}$ in Models (\protect\ref{eq:stat2-ma1})
and (\protect\ref{eq:stat2-ma2})}
\label{tab6}
\begin{tabular*}{\textwidth}{@{\extracolsep{\fill}}lccccc@{}}
\hline
&\multicolumn{5}{c@{}}{$\bolds{p}$}\\[-5pt]
&\multicolumn{5}{c@{}}{\hrulefill}\\
 $\bolds{n_1=n_2}$&\multicolumn{1}{c}{$\bolds{50}$} &
\multicolumn{1}{c}{$\bolds{100}$} &\multicolumn{1}{c}{$\bolds{200}$}
&\multicolumn{1}{c}{$\bolds{500}$}&\multicolumn{1}{c@{}}{$\bolds{700}$} \\
\hline
\multicolumn{6}{@{}c@{}}{Sizes}\\
\phantom{0}20 & 0.105 & 0.092 & 0.085 & 0.082 & 0.082\\
\phantom{0}50 & 0.101 & 0.090 & 0.081 & 0.088 & 0.090\\
\phantom{0}80 & 0.107 & 0.094 & 0.083 & 0.078 & 0.065\\
100 & 0.093 & 0.083 & 0.093 & 0.059 & 0.071\\[6pt]
\multicolumn{6}{@{}c@{}}{Powers}\\
\phantom{0}20 & 0.499 & 0.501 & 0.519 & 0.482& 0.502\\
\phantom{0}50 & 0.775 & 0.802 & 0.783 & 0.754& 0.777\\
\phantom{0}80 & 0.945 & 0.923 & 0.921 & 0.922& 0.923\\
100& 0.974 & 0.957 & 0.969 & 0.964& 0.960 \\
\hline
\end{tabular*}
\end{table}

We observed from Table \ref{tab5} that the empirical sizes of the proposed test
converged to the nominal 5\% quite rapidly, while the powers were quite
high and quickly increased to 1. For the Gamma distributed case
reported in Table \ref{tab6}, the convergence of the empirical sizes to the
nominal level was slower than the normally distributed case indicating
that the convergence of the asymptotic normality depends on the
underlying distribution, as well as the sample size and dimensionality.
The powers in Table \ref{tab6} were reasonable, although they were smaller than
the corresponding normally distributed case in Table \ref{tab5}. Nevertheless,
the power was quite responsive to the increase of~$p$ and the sample sizes.

\section{An empirical study}\label{sec5}

We report an empirical study on a leukemia data by applying the
proposed tests on\vadjust{\goodbreak} the variance--covariance matrices. The~da\-ta [Chiaretti
et al. (\citeyear{r15})], available from \url{http://www.bioconductor.org/},
consist of microarray expressions of 128 patients with either T-cell or
B-cell acute
lymphoblastic leukemia (ALL); see Dudoit, Keles and van der Laan (\citeyear{r17})
and \citet{r13} for analysis on the same dataset. We considered
a subset of the ALL data of 79 patients with the B-cell ALL. We were
interested in two types of the B-cell tumors: BCR/ABL, one of the most
frequent cytogenetic abnormalities in human leukemia, and NEG, the
cytogenetically normal B-cell ALL. The number of patients with BCR/ABL
was 37 and that with NEG was 42.

A major motivation for developing the proposed test procedures for
high-dimensional variance--covariance matrices comes from the need to
identify sets of genes which are significantly different
with respect to two treatments in genetic research; see \citet{r7},
Efron and Tibshrini (\citeyear{r19}), Newton et al. (\citeyear{r37}) and \citet{r36} for comprehensive discussions.
Biologically speaking, each gene does not function
individually, but rather tends to work with others to achieve certain
biological tasks.
Gene-sets are technically defined vocabularies which produce names of
gene-sets (also called GO terms).
There are three categories of Gene ontologies of interest: 
Biological Processes (BP), Cellular Components (CC) and Molecular
Functions (MF).
For the ALL data, a preliminary screening with gene-filtering left a
total number of 2391 genes
for analysis with 1599 unique GO terms in BP category, 290 in CC and
357 in MF.

\begin{figure}

\includegraphics{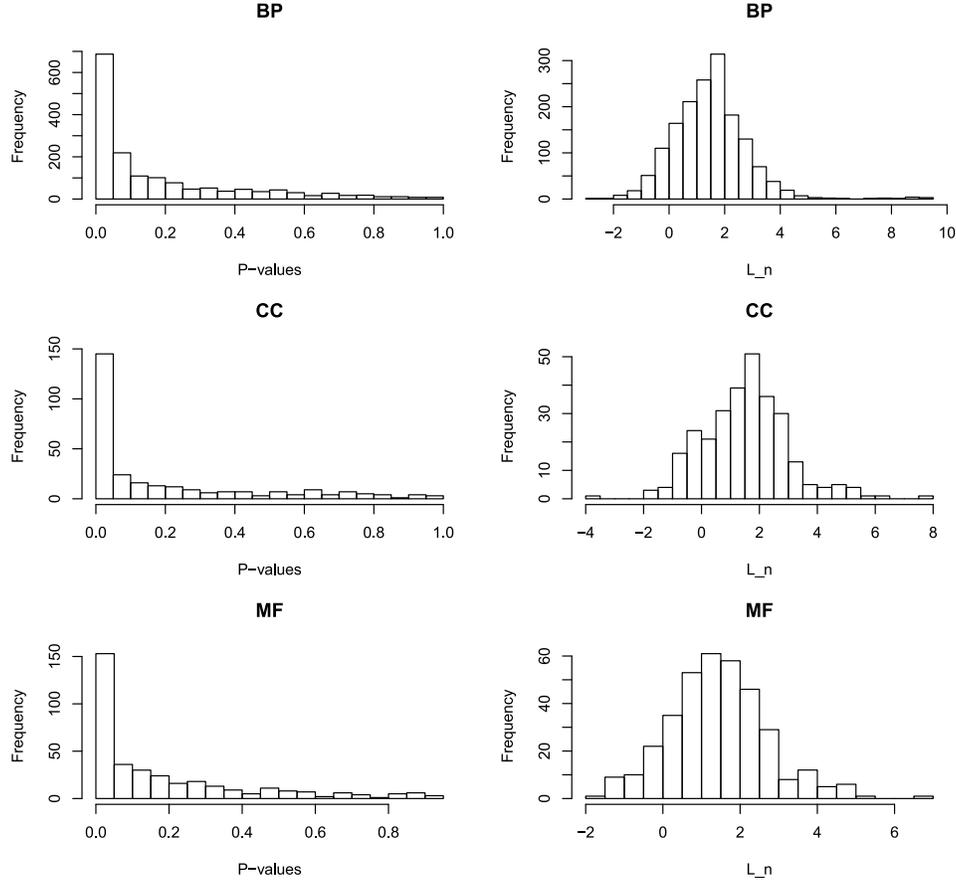}

\caption{Histograms of $p$-values (left panels) for testing two
covariance matrices and test statistic $L_n$ (right panels) for the
three gene-categories.} \label{fig1}
\end{figure}

Let us denote ${\mathcal{S}}_1, \ldots, {\mathcal{S}}_q$ for $q$ gene-sets,
where ${\mathcal{S}}_g$ consists of $p_g$ genes.
Let~$F_{1 {\mathcal{S}}_g}$ and $F_{2 {\mathcal{S}}_g}$ be the distribution functions
corresponding to ${\mathcal{S}}_g$ under the treatment and control, and
$\mu_{1 {\mathcal{S}}_g}$ and $\mu_{2 {\mathcal{S}}_g}$ be their respective
means, and $\Sigma_{1 {\mathcal{S}}_g}$ and $\Sigma_{2 {\mathcal{S}}_g}$ be
their respective variance--covariance matrices. Our first hypotheses of
interest are
$H_{0 g}\dvtx \Sigma_{1 {\mathcal{S}}_g} = \Sigma_{2 {\mathcal{S}}_g}$ {for}
$g=1, \ldots, q$ regarding the variance--covariance matrices.
For the second hypothesis, we divided each gene-set into two
sub-vectors by selecting the first $[p/2]$ dimensions of the gene-set
as the first segment and the rest as the second.

We first applied the proposed test for the equality of the entire
variance--covariance matrices and obtained the $p$-value for each
gene-set. The $p$-values and the values of the test statistics $L_n$ as
given in (\ref{null-stat}) are displayed in Figure \ref{fig1} for the three
gene-categories. By controlling the
false discovery rate [FDR, \citet{r8}] at 0.05, 338
GO terms were declared significant in the BP category, 77 in the CC
and 75 in the MF, indicating that the dependence structure among the
gene-sets was significantly different between the BCR/ABL and the NEG
ALL patients for a large number of gene sets. That a relatively large
number of gene-sets being declared significant by the proposed test
was not entirely surprising, as we observe from Figure \ref{fig1} that there
were very large number of $p$-values which were very close to 0.

For those GO terms 
which had been declared having different variance--covariance matrices,
we carried out a follow-up analysis trying to gain more details on the
differences by partitioning the variance--covariance into four blocks
in the form of (\ref{eq:block}) with $p_1=[p/2]$ and $p_2=p-p_1$.
We want to know if the difference was caused by the diagonal blocks or
the off-diagonal blocks.
The tests on the two diagonal blocks were conducted using the first
proposed test for the variance--covariance matrix but with $p_1$ or
$p_2$ dimensions, respectively. The tests on the off-diagonal blocks
were conducted by employing the second proposed test for covariances
between the two sub-vectors. The results are summarized in Table~\ref{tab7},
which provides the numbers
of gene-sets which were tested significant in the diagonal matrices
only, the off-diagonal matrix only, and both at 5\%.
There were far more gene-sets which had both diagonal and off-diagonal
matrices being significantly different, and it was less likely that the
off-diagonal matrices were different
while the diagonal matrices were otherwise. 
It was a little surprising to see that the numbers of significant
gene-sets for the diagonally-only, off-diagonal only and both in each
functional category added up to the total numbers exactly for all three
gene-categories.

\begin{table}
\caption{Number of GO terms which were tested significantly different
at the diagonal blocks, off-diagonal blocks and both diagonal and
off-diagonal blocks, respectively}
\label{tab7}
\begin{tabular*}{\textwidth}{@{\extracolsep{\fill}}lcccc@{}}
\hline
& \textbf{Diagonal only} & \textbf{Off-diagonal only} & \textbf{Both} & \textbf{Total} \\
\hline
BP & 115 & 17 & 206  & 338\\
CC & \phantom{0}26 & \phantom{0}1& \phantom{0}50  & \phantom{0}77 \\
MF & \phantom{0}22 & \phantom{0}0& \phantom{0}53  & \phantom{0}75 \\
\hline
\end{tabular*}
\end{table}

\begin{table}
\caption{Two by two classifications on the number (probability) of
GO-terms rejected/not rejected by the tests for the means and the
variances for the three functional categories, respectively}
\label{tab8}
\begin{tabular*}{\textwidth}{@{\extracolsep{\fill}}lcc@{}}
\hline
&\multicolumn{2}{c@{}}{\textbf{Mean test}}\\[-5pt]
&\multicolumn{2}{c@{}}{\hrulefill}\\
\textbf{Variance test} & \textbf{Rejected} & \textbf{Not rejected} \\
\hline
\multicolumn{3}{c}{(a) Biological Processes (BP)}\\
Rejected & \phantom{0}314 (0.196) & \phantom{0}22 (0.015) \\
Not rejected & 1000 (0.625) & 263 (0.164) \\[6pt]
\multicolumn{3}{c}{(b) Cellular Components (CC)}\\
Rejected & \phantom{00}77 (0.266) & \phantom{00}4 (0.014) \\
Not rejected & \phantom{0}164 (0.566) & \phantom{0}45 (0.154) \\[6pt]
\multicolumn{3}{c}{(c) Molecular Functions (MF)}\\
Rejected & \phantom{00}86 (0.241) & \phantom{00}1  (0.003) \\
Not rejected & \phantom{0}203 (0.568) & \phantom{0}67 (0.188) \\
\hline
\end{tabular*}
\end{table}

As we have stated in the \hyperref[sec1]{Introduction}, the proposed tests are part of
the effort in testing for high-dimensional distributions between two
treatments. However, directly testing on the distribution functions is
quite challenging due to the high dimensionality as such tests may
endure low power.
A realistic and intuitive way is to test for simpler characteristics of
the distributions, for instance testing for the means as in \citet{r5} and \citet{r13}, and the variance--covariance
as considered in this paper.
For the ALL data, in addition to testing for the variance--covariance, %
we also carried out tests for the means proposed in \citet{r13}
at a
level of 5\%.  Table \ref{tab8} contains two by two classifications on the
number and the probability of gene-sets which are
rejected/not rejected by the tests for the mean and the variance
respectively. It is observed that it is far more
likely for the means to be significantly different than the
variance--covariance, with the probability of rejection
being around $0.8$ for the means versus 0.2 to 0.3 for the covariance
for the three functional categories.
Given a gene-set which was not tested significant for the means, the
conditional probability of being tested significant for the covariance
is lower than that given a gene-set was not tested significant for the means.
These were confirmed by conducting the chi-square test for association
for the three gene-set categories, which rejected overwhelmingly (with
$p$-values all less than 0.0005) the hypothesis of no-association between
being tested significant for the mean and the variance.
For this particular dataset, the tests for the means were quite effective
in disclosing most of the differentially expressed gene-sets. However,
we do see that for Biological Processes and
Cellular Component categories, among those whose means were not
declared significantly different, there were about 10\% of gene-sets
having significant different covariance structures.
\section{Technical details}\label{sec6}

As both $T_{n_1,n_2}$ and $S_{n_1,n_2}$ are invariant under the
location transformation, we assume $\mu_i=0$ throughout this section.

\subsection{Derivations of $\operatorname{Var}(T_{n_1,n_2})$ and $\operatorname{Var}(S_{n_1,n_2})$}\label{sec6.1}

Recall that $T_{n_1,n_2}=A_{n_1}+A_{n_2}-2C_{n_1n_2}$. It is
straightforward to show that $\mathrm{E}(T_{n_1,n_2})=\operatorname{tr}\{(\Sigma
_1-\Sigma_2)^2\}$. By noticing that $\operatorname{Cov}(A_{n_1},A_{n_2})=0$,
\begin{eqnarray}\label{var-stat}
\operatorname{Var}(T_{n_1,n_2})&=&\operatorname{Var}(A_{n_1})+\operatorname
{Var}(A_{n_2})+4\operatorname{Var}(C_{n_1n_2})\nonumber\\[-8pt]\\[-8pt]
&&{}- 4\operatorname{Cov}(A_{n_1},C_{n_1n_2})- 4\operatorname
{Cov}(A_{n_2},C_{n_1n_2}).\nonumber
\end{eqnarray}

Adopting results from \citet{r14}, for $h=1$ or $2$,
\begin{eqnarray}\label{eqn7-1}
\qquad\operatorname{Var}(A_{n_h})
&=&\frac{4}{n_h^2}\operatorname{tr}^2(\Sigma_h^2)+\frac{8}{n_h}\operatorname{tr}(\Sigma
_h^4)+\frac{4\Delta_h}{n_h}
\operatorname{tr}(\Gamma_h^{\prime}\Gamma_h\Gamma_h^{\prime}\Gamma_h \circ
\Gamma_h^{\prime}\Gamma_h
\Gamma_h^{\prime}\Gamma_h)\nonumber\\[-8pt]\\[-8pt]
&&{}+O\biggl\{\frac{1}{n_h^3}\operatorname{tr}^2(\Sigma_h^2)+\frac{1}{n_h^2}\operatorname{tr}(\Sigma
_h^4)\biggr\}.\nonumber
\end{eqnarray}

Furthermore, we obtain
\begin{eqnarray}\label{var-c}
\operatorname{Var}(C_{n_1n_2})
&=&\frac{2}{n_1n_2}\operatorname{tr}^2(\Sigma_1
\Sigma_2)+\biggl(\frac{2}{n_1}+\frac{2}{n_2}\biggr)\operatorname{tr}(\Sigma_1 \Sigma_2
\Sigma_1
\Sigma_2)\hspace*{-30pt}\nonumber\\
&&{}+\frac{\Delta_1}{n_1}\operatorname{tr}(\Gamma_1^{\prime}
\Gamma_2\Gamma_2^{\prime}\Gamma_1\circ\Gamma_1^{\prime}
\Gamma_2\Gamma_2^{\prime}\Gamma_1)\nonumber\\[-8pt]\\[-8pt]
&&{}+\frac{\Delta_2}{n_2}\operatorname{tr}(\Gamma_2^{\prime}
\Gamma_1\Gamma_1^{\prime}\Gamma_2\circ\Gamma_2^{\prime}
\Gamma_1\Gamma_1^{\prime}\Gamma_2)+o\biggl\{\frac{1}{n_1n_2}\operatorname{tr}^2(\Sigma
_1 \Sigma_2)\biggr\}\nonumber\\
&&{}+ O \Biggl[\Biggl\{\frac{1}{\sqrt{n_1n_2}}+\frac{1}{n_1n_2}
+\sum_{i=1}^2\biggl(\frac{1}{\sqrt{n_i}}+\frac{1}{n_i}\biggr)\Biggr\}\operatorname
{Var}(C_{n_1n_2,1}) \Biggr].\nonumber
\end{eqnarray}

By carrying out similar procedures, we are able to obtain $\operatorname
{Cov}(A_{n_1},C_{n_1n_2})$ and $\operatorname{Cov}(A_{n_2},C_{n_1n_2})$. After
we substitute all the results into (\ref{var-stat}),
\begin{eqnarray}\label{stat1-final}
\operatorname{Var}(T_{n_1n_2})
&=&\sum_{i=1}^2 \biggl[\frac{4}{n_i^2}\operatorname{tr}^2(\Sigma_i^2)+\frac
{8}{n_i}\operatorname{tr}(\Sigma_i^4)
+\frac{4\Delta_i}{n_i}
\operatorname{tr}(\Gamma_i^{\prime}\Gamma_i\Gamma_i^{\prime}\Gamma_i \circ
\Gamma_i^{\prime}\Gamma_i\Gamma_i^{\prime}\Gamma_i)\nonumber\\
&&\hspace*{16pt}\hspace*{41pt}{}-\frac{16}{n_i}\operatorname{tr}(\Sigma_i^2\Sigma_1\Sigma_2)
-\frac{8\Delta_i}{n_i} \operatorname{tr}(\Gamma_i^{\prime}\Sigma_1\Gamma_i
\circ\Gamma_i^{\prime}\Sigma_2\Gamma_i) \biggr]\nonumber\\
&&{}+\frac{8}{n_1n_2}\operatorname{tr}^2(\Sigma_1 \Sigma_2)+\biggl(\frac{8}{n_1}+\frac
{8}{n_2}\biggr)\operatorname{tr}(\Sigma_1 \Sigma_2
\Sigma_1 \Sigma_2)\nonumber\\
&&{}+\frac{4\Delta_1}{n_1}\operatorname{tr}(\Gamma_1^{\prime}
\Gamma_2\Gamma_2^{\prime}\Gamma_1\circ\Gamma_1^{\prime}
\Gamma_2\Gamma_2^{\prime}\Gamma_1)\nonumber\\[-8pt]\\[-8pt]
&&{}+\frac{4\Delta_2}{n_2}\operatorname{tr}(\Gamma
_2^{\prime}
\Gamma_1\Gamma_1^{\prime}\Gamma_2\circ\Gamma_2^{\prime}
\Gamma_1\Gamma_1^{\prime}\Gamma_2)\nonumber\\
&&{}+o\biggl\{\frac{1}{n_1n_2}\operatorname{tr}^2(\Sigma_1
\Sigma_2)\biggr\}\nonumber\\
&&{}+O \Biggl[\Biggl\{\frac{1}{\sqrt{n_1n_2}}+\frac{1}{n_1n_2}+\sum_{i=1}^2\biggl(\frac{1}{\sqrt{n_i}}+\frac{1}{n_i}\biggr)\Biggr\}\operatorname
{Var}(C_{n_1n_2,1})\nonumber\\
&&\hspace*{112pt}{} +\sum_{i=1}^2\biggl\{\frac{1}{n_i^2}\operatorname{tr}(\Sigma
_i^4)+\frac{1}{n_i^3}\operatorname{tr}^2(\Sigma_i^2)\biggr\} \Biggr]
.\nonumber
\end{eqnarray}

Similarly to $T_{n_1,n_2}$, we have $\mathrm{E}(S_{n_1,n_2})=\operatorname{tr}\{(\Sigma
_{1,12}-\Sigma_{2,12})(\Sigma_{1,12}-\Sigma_{2,12})^{\prime} \}$
and the leading order terms in $\operatorname{Var}(S_{n_1n_2})$ are given by
%
\begin{eqnarray}\label{stat2-final}
\operatorname{Var}(S_{n_1n_2})
&=&\sum_{i=1}^2 \biggl[\frac{2}{n_i^2}\operatorname{tr}^2(\Sigma_{i,12}\Sigma
_{i,12}^{\prime})
+\frac{2}{n_i^2}\operatorname{tr}(\Sigma_{i,11}^2)\operatorname{tr}(\Sigma_{i,22}^2)\nonumber\\
&&\hspace*{17pt}{}+\frac{4}{n_i}\operatorname{tr}\{(\Sigma_{i,12}\Sigma_{1,12}^{\prime}-\Sigma
_{i,12}\Sigma_{2,12}^{\prime})^2\}\nonumber\\
&&\hspace*{17pt}{} + \frac{4}{n_i}\operatorname{tr}\{(\Sigma_{i,11}\Sigma_{1,12}-\Sigma
_{i,11}\Sigma_{2,12})
(\Sigma_{i,22}\Sigma_{1,12}^{\prime}-\Sigma_{i,22}\Sigma
_{2,12}^{\prime})\}\nonumber\\[-8pt]\\[-8pt]
&&\hspace*{17pt}{} + \frac{4\Delta_i}{n_i}\nonumber\\
&&\hspace*{29pt}{}\times\operatorname{tr}\bigl\{{\Gamma_i^{(1)\prime}}(\Sigma
_{1,12}-\Sigma_{2,12})\Gamma_i^{(2)}
\circ{\Gamma_i^{(1)\prime}}(\Sigma_{1,12}-\Sigma_{2,12})\Gamma
_i^{(2)}\bigr\} \biggr]\nonumber\\
&&{}+\frac{4}{n_1n_2}
\operatorname{tr}^2(\Sigma_{1,12}\Sigma_{2,12}^{\prime})
+ \frac{4}{n_1n_2}
\operatorname{tr}(\Sigma_{1,11}\Sigma_{2,11})\operatorname{tr}(\Sigma_{1,22}\Sigma_{2,22}).\nonumber
\end{eqnarray}

\subsection{\texorpdfstring{Proof of Theorem \protect\ref{teo1}}{Proof of Theorem 1}}\label{sec6.2}
The leading order terms in $\operatorname{Var}(T_{n_1,n_2})$ are contributed by
$A_{n_h,1}$ for $h=1,2$ and $C_{n_1n_2,1}$, which are defined by
\begin{eqnarray*}
A_{n_h,1}=\frac{1}{n_h(n_h-1)} \sum_{i \ne j} (X_{hi}^{\prime}
X_{hj})^2,\qquad
C_{n_1n_2,1}
=
\frac{1}{n_1n_2}\sum_{ij}(X_{1i}^{\prime}X_{2j})^2.
\end{eqnarray*}
Hence, we only need to study the asymptotic normality of $Z_{n_1,n_2}$
which is defined by $
Z_{n_1,n_2}=: A_{n_1,1}+A_{n_2,1}-2C_{n_1n_2,1}$.

In order to construct a martingale
sequence, it is convenient to have new random variables
${Y_i}$ which are defined as
\begin{eqnarray*}
Y_i&=& X_{1i} \qquad \mbox{for }   i=1,2,\ldots,n_1,\\
Y_{n_1+j}&=& X_{2j}\qquad  \mbox{for }   j=1,2,\ldots,n_2.
\end{eqnarray*}

To construct a martingale difference, we let
$\mathscr{F}_0=\{\varnothing, \Omega\}$, $\mathscr{F}_k=\sigma\{Y_1,\ldots,\allowbreak Y_k
\}$ with $k=1,2,\ldots,n_1+n_2$. And let $\mathrm{E}_k(\cdot)$ denote the
conditional expectation given $\mathscr{F}_k$. Define
$D_{n,k}=(\mathrm{E}_k-\mathrm{E}_{k-1})Z_{n_1,n_2}$ and it is easy to see
that
$Z_{n_1,n_2}-\mathrm{E}(Z_{n_1,n_2})=\sum_{k=1}^{n_1+n_2}D_{n,k}$.

\begin{lemma}\label{lem1}
 For any $n$, $\{D_{n,k}, 1\le k \le n\}$ is a
martingale difference sequence with respect to the $\sigma$-fields
$\{\mathscr{F}_k,1\le k \le n\}$.
\end{lemma}

\begin{pf}
 First of all, it is straightforward to show that $\mathrm{E}D_{n,k}=0$. Next, by denoting $S_{n,m}=
\sum_{k=1}^{m}D_{n,k}=\mathrm{E}_mZ_{n_1,n_2}-\mathrm{E}Z_{n_1,n_2}$, we have
$S_{n,q}=S_{n,m}+(\mathrm{E}_qZ_{n_1,n_2}
-\mathrm{E}_mZ_{n_1,n_2})$. Then we can show that $\mathrm{E}(S_{n,q}|\mathscr{F}_m)=S_{n,m}$.
This completes the proof of Lemma \ref{lem1}.
\end{pf}

To apply martingale central limit theorem, we need Lemmas \ref{lem2} and \ref{lem3}.

\begin{lemma}\label{lem2}
Under Condition \textup{A2} and as $\min \{n_1,n_2\}\to
\infty$,
\begin{eqnarray*}
\frac{\sum_{k=1}^{n_1+n_2}\sigma^2_{n,k}}{\operatorname
{Var}(Z_{n_1,n_2})}\stackrel{p}{\rightarrow}  1,
\end{eqnarray*}
where $\sigma^2_{n,k}=\mathrm{E}_{k-1}(D_{n,k}^2)$.
\end{lemma}

\begin{pf}
To prove Lemma \ref{lem2}, first we can show $\mathrm{E}(\sum
_{k=1}^{n_1+n_2}\sigma_{n,k}^2)=\break\operatorname{Var}(Z_{n_1,n_2})$. Then we
will show that as $\min \{n_1,n_2\}\to\infty$, $\operatorname
{Var}(\sum_{k=1}^{n_1+n_2}\sigma_{n,k}^2)/\allowbreak\operatorname
{Var}^2(Z_{n_1,n_2})\to
0$. To this end, we decompose $\sum_{k=1}^{n_1+n_2}\sigma_{n,k}^2$
into the sum of eight parts,
\begin{eqnarray*}
\sum_{k=1}^{n_1+n_2}\sigma
_{n,k}^2=R_1+R_2+R_3+R_4+R_5+R_6+R_7+R_8,
\end{eqnarray*}
where with $Q_{1,k-1}=\sum_{i=1}^{k-1}(Y_iY_i^{\prime}-\Sigma_1)$
and $Q_{2,n_1+l-1}=\sum_{i=1}^{l-1}(Y_{n_1+i}Y_{n_1+i}^{\prime
}-\Sigma_2)$,
\begin{eqnarray*}
R_1&=&
\sum_{k=1}^{n_1}\frac{8}{n_1^2(n_1-1)^2}\operatorname{tr}(Q_{1,k-1}\Sigma
_1Q_{1,k-1}\Sigma_1)\\
&&{}+\sum_{l=1}^{n_2}\frac{8}{n_2^2(n_2-1)^2}\operatorname{tr}(Q_{2,n_1+l-1}\Sigma
_2Q_{2,n_1+l-1}\Sigma_2),\\
R_2&=&\sum_{k=1}^{n_1}\frac{16}{n_1^2(n_1-1)}\sum_{i=1}^{k-1}\{
Y_i^{\prime}(\Sigma_1^3-\Sigma_1
\Sigma_2\Sigma_1)Y_i \},\\
R_3&=&\sum_{l=1}^{n_2}\frac{16}{n_2^2(n_2-1)}
\Biggl[\operatorname{tr}(Q_{2,n_1+l-1}\Sigma_2^3)-\operatorname{tr}\Biggl\{Q_{2,n_1+l-1}\Sigma_2
\Biggl(\frac{1}{n_1}\sum_{i=1}^{n_1}Y_iY_i^{\prime}\Biggr)\Sigma_2\Biggr\}
\Biggr],\\
R_4&=&\frac{8}{n_1^2
n_2}\sum_{i,j}^{n_1}\operatorname{tr}(Y_jY_j^{\prime}\Sigma_2Y_iY_i^{\prime}\Sigma
_2)-\frac{16}{n_1
n_2}\operatorname{tr}\Biggl\{\Sigma_2^3\Biggl(\sum_{i=1}^{n_1}Y_iY_i^{\prime}\Biggr)\Biggr\},\\
R_5&=&\sum_{k=1}^{n_1}\frac{4\Delta_1}{n_1^2(n_1-1)^2}\operatorname{tr}(\Gamma
_1^{\prime}Q_{1,k-1}\Gamma_1
\circ
\Gamma_1^{\prime}Q_{1,k-1}\Gamma_1)\nonumber\\
&&{}+\sum_{l=1}^{n_2}\frac{4\Delta_2}{n_2^2(n_2-1)^2}\operatorname{tr}(\Gamma
_2^{\prime}Q_{2,n_1+l-1}\Gamma_2
\circ\Gamma_2^{\prime}Q_{2,n_1+l-1}\Gamma_2),\\
R_6&=&\sum_{k=1}^{n_1}\frac{8\Delta_1}{n_1^2(n_1-1)}\operatorname{tr}\{\Gamma
_1^{\prime}(\Sigma_1-\Sigma_2)\Gamma_1
\circ\Gamma_1^{\prime}Q_{1,k-1}\Gamma_1\},\\
R_7&=&\sum_{l=1}^{n_2}\frac{8\Delta_2}{n_2^2(n_2-1)}
\Biggl[\operatorname{tr}(\Gamma_2^{\prime}Q_{2,n_1+l-1}\Gamma_2
\circ\Gamma_2^{\prime}\Sigma_2
\Gamma_2)\\[-4pt]
&&\hspace*{70pt}{}-\operatorname{tr}\Biggl\{\Gamma_2^{\prime}Q_{2,n_1+l-1}\Gamma_2\circ
\Gamma_2^{\prime}\Biggl(\frac{1}{n_1}\sum_{i=1}^{n_1}Y_iY_i^{\prime
}\Biggr)\Gamma_2\Biggr\} \Biggr]
\end{eqnarray*}
and
\begin{eqnarray*}
R_8&=&\frac{4\Delta_2}{n_1^2
n_2}\sum_{i,j}^{n_1}\operatorname{tr}(\Gamma_2^{\prime}Y_iY_i^{\prime}\Gamma
_2\circ
\Gamma_2^{\prime}Y_jY_j^{\prime}\Gamma_2)-\frac{8\Delta_2}{n_1
n_2}\sum_{i=1}^{n_1}\operatorname{tr}(\Gamma_2^{\prime}\Sigma_2\Gamma_2\circ
\Gamma_2^{\prime}Y_iY_i^{\prime}\Gamma_2).
\end{eqnarray*}

Therefore, we need to show that
$\operatorname{Var}(R_i)=o\{\operatorname{Var}^2(Z_{n_1,n_2})\}$ for $i=1,\ldots,8$.

For $R_1$, there exists a constant $K_1$ such that
\begin{eqnarray*}
\operatorname{Var}(R_1)\le K_{1}\{
n_1^{-4}\operatorname{tr}^2(\Sigma_1^2)\operatorname{tr}(\Sigma_1^4)+
n_2^{-4}\operatorname{tr}^2(\Sigma_2^2)\operatorname{tr}(\Sigma_2^4)\}.
\end{eqnarray*}

Then, applying $\operatorname{Var}^2(Z_{n_1,n_2})\ge
\frac{16}{n_1^4}\operatorname{tr}^4(\Sigma_1^2)+\frac{16}{n_2^4}\operatorname{tr}^4(\Sigma_2^2)$
from (\ref{variance}), we know
\begin{eqnarray*}
\frac{\operatorname{Var}(R_1)}{\operatorname{Var}^2(Z_{n_1,n_2})}
\le\frac{K_1}{16}\biggl \{
\frac{\operatorname{tr}(\Sigma_1^4)}{\operatorname{tr}^2(\Sigma_1^2)}+
\frac{\operatorname{tr}(\Sigma_2^4)}{\operatorname{tr}^2(\Sigma_2^2)} \biggr\},
\end{eqnarray*}
where $\operatorname{tr}(\Sigma_1^4)/\operatorname{tr}^2(\Sigma_1^2)\! \to\!0$ under Condition A2.
Thus, $\operatorname{Var}(R_1)=\break o\{\operatorname{Var}^2(Z_{n_1,n_2})\}$.

By carrying out similar procedures we can show that the above is
true for~$R_i$ with $i=1,\ldots,8$. Hence we complete the proof of Lemma \ref{lem2}.
\end{pf}

\begin{lemma}\label{lem3}
Under Condition \textup{A2}, as $\min \{n_1,n_2\}\to
\infty$
\begin{eqnarray*}
\frac{\sum_{k=1}^{n_1+n_2}\mathrm{E}(D_{n,k}^4)}{\operatorname
{Var}^2(Z_{n_1,n_2})}\to
0.
\end{eqnarray*}
\end{lemma}

\begin{pf}
For the case of $1\le k\le n_1$, there exists a constant $c$
such that
\begin{eqnarray*}
\sum_{k=1}^{n_1}\mathrm{E}(D_{n,k}^4)
\le
c [n_1^{-3}\operatorname{tr}^2\{(\Sigma_1^2-\Sigma_1\Sigma_2)^2\}
+n_1^{-5}\operatorname{tr}^4\{(\Sigma_1^2)\}
 ].
\end{eqnarray*}

Using the results $\operatorname{Var}^2(Z_{n_1,n_2})\ge64 n_1^{-2}\operatorname{tr}^2\{
(\Sigma_1^2-\Sigma_1\Sigma_2)^2\}$ and\break $\operatorname
{Var}^2(Z_{n_1,n_2})\ge16n_1^{-4}\operatorname{tr}^4\{(\Sigma_1^2)\}$ from (\ref{variance})
and as $n_1 \to\infty$, we have
\begin{eqnarray*}
\frac{\sum_{k=1}^{n_1}\mathrm{E}(D_{n,k}^4)}{\operatorname
{Var}^2(Z_{n_1,n_2})} \le\frac{c}{n_1} \to0.
\end{eqnarray*}

For the case of $n_1<k <n_1+n_2$, there exists a constant $d$ such
that
\begin{eqnarray}\label{result2}\qquad
\sum_{k=n_1}^{n_1+n_2}\mathrm{E}(D_{n,k}^4)
&\le&\frac{d}{n_1^2n_2^4}\{2\operatorname{tr}^4(\Sigma_1\Sigma_2)+\operatorname{tr}^2(\Sigma
_1\Sigma_2)\operatorname{tr}^2(\Sigma_1^2)\}\nonumber\\
&&{}+\frac{d}{n_1n_2^4} [2\operatorname{tr}^2(\Sigma_1\Sigma_2)\operatorname{tr}\{(\Sigma
_2^2-\Sigma_2\Sigma_1)^2 \} ]+\frac{d}{n_2^5}\operatorname{tr}^4\{(\Sigma
_2^2)\}
\\[-2pt]
&&{}+\frac{d}{n_2^4} [2\operatorname{tr}^2(\Sigma_2^2)\operatorname{tr}\{(\Sigma_2^2-\Sigma
_2\Sigma_1)^2 \}+4\operatorname{tr}^2(\Sigma_1\Sigma_2)\operatorname{tr}^2(\Sigma_2^2)
].\nonumber
\end{eqnarray}

To evaluate the ratio of individual term in (\ref{result2}) to $\operatorname
{Var}^2(Z_{n_1,n_2})$, respectively, we simply replace $\operatorname
{Var}^2(Z_{n_1,n_2})$ by corresponding terms in (\ref{variance}). Then
under Condition A2 and as $n_2 \to\infty$,
$\sum_{k=n_1+1}^{n_1+n_2}\mathrm{E}(D_{n,k}^4)/\operatorname
{Var}^2(Z_{n_1,n_2})\to0$. Therefore, we complete the proof of Lemma \ref{lem3}.
\end{pf}

With two sufficient conditions given in Lemmas \ref{lem2} and \ref{lem3}, we
conclude that
\begin{eqnarray*}
\frac{Z_{n_1,n_2}-\mathrm{E}(Z_{n_1,n_2})}{\operatorname{Var}(Z_{n_1,n_2})}
\stackrel{d}{\rightarrow}  \mathrm{N}(0,1).
\end{eqnarray*}

If we let $\varepsilon
_{n_1,n_2}=A_{n_1,2}+A_{n_1,3}+A_{n_2,2}+A_{n_2,3}-2C_{n_1n_1,2}-2C_{n_1n_1,3}-2C_{n_1n_1,4}$,
then $T_{n_1,n_2}=Z_{n_1,n_2}+\varepsilon_{n_1,n_2}$. From $\operatorname
{Var}(\varepsilon_{n_1,n_2})=o(\sigma_{n_1,n_2}^2)$,
\begin{eqnarray*}
\operatorname{Var}\biggl(\frac{\varepsilon_{n_1,n_2}}{\sigma_{n_1,n_2}}\biggr)&=&\frac
{\operatorname{Var}(\varepsilon_{n_1,n_2})}{\sigma_{n_1,n_2}^2}
\to0.
\end{eqnarray*}

Moreover, $\mathrm{E}(\varepsilon_{n_1,n_2})=0$. Therefore, $\varepsilon
_{n_1,n_2}/\sigma_{n_1,n_2}\stackrel{p}{\rightarrow}  0$. From Slutsky's
Theorem, we complete the proof of Theorem \ref{teo1}.

\subsection{\texorpdfstring{Proof of Theorem \protect\ref{teo2}}{Proof of Theorem 2}}\label{sec6.3}
Recall that $\mathrm{E}(A_{n_h})=\operatorname{tr}(\Sigma_h^2)$ for $h=1$ or $2$. To
show $A_{n_h}/\operatorname{tr}(\Sigma_h^2)\stackrel{p}{\rightarrow}  1$, it is sufficient to
show that $\operatorname{Var}\{A_{n_h}/\operatorname{tr}(\Sigma_h^2)\}\to0$.

From (\ref{eqn7-1}), we have
\begin{eqnarray*}
& &\operatorname{Var} \biggl\{\frac{A_{n_h}}{\operatorname{tr}(\Sigma_h^2)} \biggr\}
\\
&&\qquad{}\le \frac{1}{\operatorname{tr}^2(\Sigma_h^2)} \biggl[\frac{4}{n_h^2}\operatorname{tr}^2(\Sigma
_h^2)+\frac{8+4\Delta_h}{n_h}\operatorname{tr}(\Sigma_h^4)
+O\biggl\{\frac{1}{n_h^3}\operatorname{tr}^2(\Sigma_h^2)+\frac{1}{n_h^2}\operatorname{tr}(\Sigma_h^4)\biggr\}
 \biggr],
\end{eqnarray*}
where $\operatorname{tr}(\Sigma_h^4)/\operatorname{tr}^2(\Sigma_h^2)\to0$ under Condition A2.
Hence, $A_{n_h}/\operatorname{tr}(\Sigma_h^2)\stackrel{p}{\rightarrow}  1$.

Moreover, under $H_{0 a}\dvtx \Sigma_1=\Sigma_2=\Sigma$,
$A_{n_h}/\operatorname{tr}(\Sigma^2)\stackrel{p}{\rightarrow}  1$. Then using the continuous
mapping theorem, we have $\hat{\sigma}_{0,n_1,n_2}/\sigma
_{0,n_1,n_2}\stackrel{p}{\rightarrow}  1$.

\subsection{\texorpdfstring{Proof of Theorem \protect\ref{teo3}}{Proof of Theorem 3}}\label{sec6.4}
The leading order terms in $\operatorname{Var}(S_{n_1,n_2})$ are contributed
by $U_{n_h,1}$ and $W_{n_1n_2,1}$ which are defined by
\begin{eqnarray*}
U_{n_h,1}&=&\frac{1}{n_h(n_h-1)} \sum_{i \ne j}
{X_{hi}^{(1)\prime}}X_{hj}^{(1)}
{X_{hj}^{(2)\prime}}X_{hi}^{(2)},\\
W_{n_1n_2,1}&=&
\frac{1}{n_1n_2}\sum_{ij}{X_{1i}^{(1)\prime}}
X_{2j}^{(1)}{X_{2j}^{(2)\prime}}X_{1i}^{(2)}.
\end{eqnarray*}
From Slutsky's Theorem, we only need to study the asymptotic normality
of $H_{n_1,n_2}$ which is defined as $H_{n_1,n_2}=:
U_{n_1,1}+U_{n_2,1}-2W_{n_1n_2,1}$.

To implement martingale central limit theorem to $H_{n_1,n_2}$, we need
a martingale
sequence. To this end, we define random variables which are
\begin{eqnarray*}
Y_i^{(1)}&=& X_{1i}^{(1)} \quad \mbox{and} \quad  Y_i^{(2)}=
X_{1i}^{(2)}
\qquad\mbox{for }   i=1,2,\ldots,n_1,\\
Y_{n_1+j}^{(1)}&=& X_{2j}^{(1)} \quad \mbox{and}\quad
Y_{n_1+j}^{(2)}= X_{2j}^{(2)} \qquad\mbox{for }
j=1,2,\ldots,n_2.
\end{eqnarray*}

If we define $C_{n,k}=(\mathrm{E}_k-\mathrm{E}_{k-1})H_{n_1,n_2}$, where
$\mathrm{E}_k(\cdot)$ denote the conditional expectation given
$\mathscr{F}_k=\sigma\{Y_1,\ldots,Y_k
\}$ with $k=1,2,\ldots,n_1+n_2$, we claim that $\{C_{n,k}, 1\le k \le n\}$
is a martingale difference sequence with respect to the $\sigma$-fields
$\{\mathscr{F}_k,1\le k \le n\}$ from Lemma \ref{lem1}. We need Lemmas \ref{lem4} and \ref{lem5}
to implement the martingale central limit theorem.

\begin{lemma}\label{lem4}
Under Conditions \textup{A2} and \textup{A4}, as $\min \{
n_1,n_2\}\to\infty$,
\begin{eqnarray*}
\frac{\sum_{k=1}^{n_1+n_2}\tau^2_{n,k}}{\operatorname
{Var}(H_{n_1,n_2})}\stackrel{p}{\rightarrow}  1,
\end{eqnarray*}
where $\tau^2_{n,k}=\mathrm{E}_{k-1}(C_{n,k}^2)$.
\end{lemma}

\begin{pf}
First, we can show that $\mathrm{E}(\sum_{k=1}^{n_1+n_2}\tau
^2_{n,k})=\operatorname{Var}(H_{n_1,n_2})$. Therefore, we only need to show
$\operatorname{Var}(\sum_{k=1}^{n_1+n_2}\tau^2_{n,k})=o\{\operatorname
{Var}^2(H_{n_1,n_2})\}$ to complete the proof of Lemma \ref{lem4}. To this end,
we decompose $\sum_{k=1}^{n_1+n_2}\tau_{n,k}^2$ into twelve parts,
\begin{eqnarray*}
\sum_{k=1}^{n_1+n_2}\sigma_{n,k}^2
=P_1+P_2+P_3+P_4+P_5+P_6
+P_7+P_8+P_9+P_{10}+P_{11}+P_{12},
\end{eqnarray*}
where with
\begin{eqnarray*}
O_{1,k-1}&=&\sum_{i=1}^{k-1}\bigl(Y_i^{(1)}{Y_i^{(2)\prime}}-\Sigma
_{1,12}\bigr) \quad\mbox{and}\\
   O_{2,n_1+l-1}&=&\sum
_{i=1}^{l-1}\bigl(Y_{n_1+i}^{(1)}{Y_{n_1+i}^{(2)\prime}}-\Sigma_{2,12}\bigr),
\\
P_1&=&
\sum_{k=1}^{n_1}\frac{4}{n_1^2(n_1-1)^2}\operatorname{tr}(O_{1,k-1}\Sigma
_{1,12}^{\prime}
O_{1,k-1}\Sigma_{1,12}^{\prime})\\
&&{}+\sum_{l=1}^{n_2}\frac{4}{n_2^2(n_2-1)^2}\operatorname{tr}(O_{2,n_1+l-1}\Sigma
_{2,12}^{\prime}
O_{2,n_1+l-1}\Sigma_{2,12}^{\prime}),\\
P_2&=&
\sum_{k=1}^{n_1}\frac{4}{n_1^2(n_1-1)^2}\operatorname{tr}(O_{1,k-1}\Sigma_{1,22}
O_{1,k-1}^{\prime}\Sigma_{1,11})\\
&&{}+\sum_{l=1}^{n_2}\frac{4}{n_2^2(n_2-1)^2}\operatorname{tr}(O_{2,n_1+l-1}\Sigma_{2,22}
O_{2,n_1+l-1}^{\prime}\Sigma_{2,11}),\\
P_3&=&\sum_{k=1}^{n_1}\frac{8}{n_1^2(n_1-1)}\operatorname{tr}\{O_{1,k-1}\Sigma
_{1,12}^{\prime} (\Sigma_{1,12}-\Sigma_{2,12})\Sigma_{1,12}^{\prime
}\},\\
P_4&=&\sum_{k=1}^{n_1}\frac{8}{n_1^2(n_1-1)}\operatorname{tr}\{O_{1,k-1}\Sigma
_{1,22} (\Sigma_{1,12}^{\prime}-\Sigma_{2,12}^{\prime})\Sigma
_{1,11}\},\\
P_5&=&\sum_{l=1}^{n_2}\frac{8}{n_2^2(n_2-1)}\operatorname{tr}\Biggl\{O_{2,n_1+l-1}\Sigma
_{2,12}^{\prime} \Biggl(\Sigma_{2,12}\!-\!\frac{1}{n_1}\sum_{i=1}^{n_1}Y_i^{(1)}
{Y_i^{(2)\prime}}\Biggr)\Sigma_{2,12}^{\prime}\Biggr\},\\
P_6&=&\sum_{l=1}^{n_2}\frac{8}{n_2^2(n_2-1)}\operatorname{tr}\Biggl\{O_{2,n_1+l-1}\Sigma
_{2,22} \Biggl(\Sigma_{2,12}^{\prime}\!-\!\frac{1}{n_1}\sum_{i=1}^{n_1}Y_i^{(2)}
{Y_i^{(1)\prime}}\Biggr)\Sigma_{2,11}\Biggr\},\\
P_7&=&\frac{4}{n_2}\operatorname{tr}\Biggl\{\Biggl(\Sigma_{2,12}-\frac{1}{n_1}\sum_{i=1}^{n_1}Y_i^{(1)}
{Y_i^{(2)\prime}}\Biggr)\Sigma_{2,12}^{\prime}\\
&&\hspace*{28pt}{}\times\Biggl(\Sigma_{2,12}-\frac
{1}{n_1}\sum_{i=1}^{n_1}Y_i^{(1)}
{Y_i^{(2)\prime}}\Biggr)\Sigma_{2,12}^{\prime} \Biggr\},\\
P_8&=&\frac{4}{n_2}\operatorname{tr}\Biggl\{\Biggl(\Sigma_{2,12}-\frac{1}{n_1}\sum_{i=1}^{n_1}Y_i^{(1)}
{Y_i^{(2)\prime}}\Biggr)\Sigma_{2,22}\\
&&\hspace*{28pt}{}\times\Biggl(\Sigma_{2,12}^{\prime}-\frac
{1}{n_1}\sum_{i=1}^{n_1}Y_i^{(2)}
{Y_i^{(1)\prime}}\Biggr)\Sigma_{2,11}\Biggr\},
\\
P_9&=&\sum_{k=1}^{n_1}\frac{4\Delta_1}{n_1^2(n_1-1)^2}\operatorname{tr}\bigl({\Gamma
_1^{(1)\prime}}
O_{1,k-1}\Gamma_1^{(2)}
\circ
{\Gamma_1^{(1)\prime}}
O_{1,k-1}\Gamma_1^{(2)}\bigr)\\
&&{}+\sum_{l=1}^{n_2}\frac{4\Delta_2}{n_2^2(n_2-1)^2}\operatorname{tr}\bigl({\Gamma
_2^{(1)\prime}}
O_{2,n_1+l-1}\Gamma_2^{(2)}
\circ{\Gamma_2^{(1)\prime}}
O_{2,n_1+l-1}\Gamma_2^{(2)}\bigr),\\
P_{10}&=&\sum_{k=1}^{n_1}\frac{8\Delta_1}{n_1^2(n_1-1)}\operatorname{tr}\bigl\{{\Gamma
_1^{(1)\prime}}
(\Sigma_{1,12}-\Sigma_{2,12})\Gamma_1^{(2)}
\circ{\Gamma_1^{(1)\prime}}O_{1,k-1}\Gamma_1^{(2)}\bigr\},\\
P_{11}&=&\sum_{l=1}^{n_2}\frac{8\Delta_2}{n_2^2(n_2-1)}\\
&&\hspace*{16pt}{}\times\operatorname{tr}\Biggl\{{\Gamma
_2^{(1)\prime}}
\Biggl(\Sigma_{2,12}-\sum_{i=1}^{n_1}\frac{Y_i^{(1)}
{Y_i^{(2)\prime}}}{n_1}\Biggr)\Gamma_2^{(2)} \circ{\Gamma
_2^{(1)\prime}}O_{2,n_1+l-1}\Gamma_2^{(2)}\Biggr\},
\\
P_{12}&=&\frac{4\Delta_2}{n_2}\operatorname{tr}\Biggl\{{\Gamma_2^{(1)\prime}}
\Biggl(\Sigma_{2,12}-\sum_{i=1}^{n_1}\frac{Y_i^{(1)}
{Y_i^{(2)\prime}}}{n_1}\Biggr)\Gamma_2^{(2)}\\
&&\hspace*{36pt} {}\circ{\Gamma
_2^{(1)\prime}}
\Biggl(\Sigma_{2,12}-\sum_{i=1}^{n_1}\frac{Y_i^{(1)}
{Y_i^{(2)\prime}}}{n_1}\Biggr)\Gamma_2^{(2)}\Biggr\}.
\end{eqnarray*}

For $P_1$, there exists a constant $J_1$ such that
\begin{eqnarray*}
\operatorname{Var}(P_{1})
&\le&\sum_{h=1}^2 \frac{J_1}{n_h^{4}}\{\operatorname{tr}^2(\Sigma_{h,12}\Sigma
_{h,12}^{\prime})
\operatorname{tr}(\Sigma_{h,11}\Sigma_{h,12}\Sigma_{h,22}\Sigma_{h,12}^{\prime
})\\
&&\hspace*{31pt}{}+\operatorname{tr}(\Sigma_{h,11}^2)
\operatorname{tr}(\Sigma_{h,22}^2)\operatorname{tr}(\Sigma_{h,11}\Sigma_{h,12}\Sigma_{h,22}\Sigma
_{h,12}^{\prime})\\
&&\hspace*{78pt}\hspace*{31pt}{}+\operatorname{tr}^2(\Sigma_{h,11}\Sigma_{h,12}\Sigma_{h,22}\Sigma
_{h,12}^{\prime})\}.
\end{eqnarray*}

Using $\operatorname{Var}^2(H_{n_1,n_2})\ge\frac{8}{n_h^4}\operatorname{tr}(\Sigma
_{h,11}^2)\operatorname{tr}(\Sigma_{h,22}^2)\operatorname{tr}^2(\Sigma_{h,12}\Sigma_{h,12}^{\prime
})$ from (\ref{stat2-variance}),
\begin{eqnarray*}
&&\frac{({J_1}/{(n_h^4)})\operatorname{tr}^2(\Sigma_{h,12}\Sigma_{h,12}^{\prime})
\operatorname{tr}(\Sigma_{h,11}\Sigma_{h,12}\Sigma_{h,22}\Sigma_{h,12}^{\prime
})}{\operatorname{Var}^2(H_{n_1,n_2})}\\
&&\qquad\le\frac{J_1 \operatorname{tr}(\Sigma_{h,11}\Sigma_{h,12}\Sigma_{h,22}\Sigma
_{h,12}^{\prime})}{8\operatorname{tr}(\Sigma_{h,11}^2)
\operatorname{tr}(\Sigma_{h,22}^2)},
\end{eqnarray*}
which goes to zero under Condition A4 for $h=1$ or $2$.

Similarly, using $\operatorname{Var}^2(H_{n_1,n_2})\ge
\frac{4}{n_h^4}\operatorname{tr}^2(\Sigma_{h,11}^2)\operatorname{tr}^2(\Sigma_{h,22}^2)$ from
(\ref{stat2-variance}),
\begin{eqnarray*}
\frac{J_1}{n_h^4}\operatorname{tr}^2(\Sigma_{h,11}\Sigma_{h,12}\Sigma_{h,22}\Sigma
_{h,12}^{\prime})
/\operatorname{Var}^2(H_{n_1,n_2})&\to&0,\quad \mbox{and}\\
\frac{J_1}{n_h^4}\operatorname{tr}(\Sigma_{h,11}^2)
\operatorname{tr}(\Sigma_{h,22}^2)\operatorname{tr}(\Sigma_{h,11}\Sigma_{h,12}\Sigma_{h,22}\Sigma
_{h,12}^{\prime})/
\operatorname{Var}^2(H_{n_1,n_2})&\to&0.
\end{eqnarray*}
Hence, $\operatorname{Var}(P_{1})=o\{\operatorname{Var}^2(H_{n_1,n_2})\}$. Similarly,
we have $\operatorname{Var}(P_{i})=\break o\{\operatorname{Var}^2(H_{n_1,n_2})\}$ for
$i=1,\ldots,12$. Therefore, we complete the proof of Lem\-ma~\ref{lem4}.
\end{pf}

\begin{lemma}\label{lem5}
Under Conditions \textup{A2} and \textup{A4}, as $\min \{
n_1,n_2\}\to\infty$
\begin{eqnarray*}
\frac{\sum_{k=1}^{n_1+n_2}\mathrm{E}(C_{n,k}^4)}{\operatorname
{Var}^2(H_{n_1,n_2})}\to
0.
\end{eqnarray*}
\end{lemma}

\begin{pf}
For the case of $1\le k\le n_1$, there exists a constant $c$
such that
\begin{eqnarray*}
\sum_{k=1}^{n_1}\mathrm{E}(C_{n,k}^4)
&\le&
c [n_1^{-3}\operatorname{tr}^2\{\Sigma_{1,11}(\Sigma_{1,12}-\Sigma
_{2,12})\Sigma_{1,22}(\Sigma_{1,12}^{\prime}
-\Sigma_{2,12}^{\prime})\}\\
&&\hspace*{108pt}\hspace*{6pt}{}+n_1^{-5}\operatorname{tr}^2(\Sigma_{1,11}^2)\operatorname{tr}^2(\Sigma_{1,22}^2)
 ].
\end{eqnarray*}

Applying $\operatorname{Var}^2(H_{n_1,n_2})\ge16 n_1^{-2}\operatorname{tr}^2\{\Sigma
_{1,11}(\Sigma_{1,12}-\Sigma_{2,12})\Sigma_{1,22}(\Sigma
_{1,12}^{\prime}
-\Sigma_{2,12}^{\prime})\}$ and $\operatorname{Var}^2(H_{n_1,n_2})\ge
4n_1^{-4}\operatorname{tr}^2(\Sigma_{1,11}^2)\operatorname{tr}^2(\Sigma_{1,22}^2)$ from (\ref{stat2-variance}) and as $n_1 \to\infty$,
\begin{eqnarray*}
\frac{\sum_{k=1}^{n_1}\mathrm{E}(C_{n,k}^4)}{\operatorname
{Var}^2(H_{n_1,n_2})} \le\frac{c}{n_1} \to0.
\end{eqnarray*}

For the case of $n_1 < k \le n_1+n_2$, we can find a constant $d$ such that
\begin{eqnarray}\label{stat2-result2}
& &\sum_{k=n_1}^{n_1+n_2}\mathrm{E}(C_{n,k}^4)\nonumber\\[-2pt]
&&\qquad\le\frac{d}{n_1^3n_2^3}\operatorname{tr}(\Sigma_{1,11}\Sigma_{2,11})\operatorname{tr}(\Sigma
_{1,22}\Sigma_{2,22})
\operatorname{tr}(\Sigma_{2,11}^2)\operatorname{tr}(\Sigma_{2,22}^2)\nonumber\\[-2pt]
&&\quad\qquad{}+\frac{d}{n_2^3}\operatorname{tr}^2\{(\Sigma_{2,11}\Sigma_{2,12}-\Sigma_{2,11}
\Sigma_{1,12})(\Sigma_{2,22}\Sigma_{2,12}^{\prime}-\Sigma
_{2,22}\Sigma_{1,12}^{\prime})\}\nonumber\\[-9pt]\\[-9pt]
&&\quad\qquad{}+\frac{d}{n_1n_2^3}\operatorname{tr}(\Sigma_{1,11}\Sigma_{2,11})\operatorname{tr}(\Sigma
_{1,22}\Sigma_{2,22})\nonumber\\[-2pt]
&&\hspace*{46pt}{} \times \operatorname{tr}\{\Sigma_{2,11}(\Sigma_{2,12}-
\Sigma_{1,12})\Sigma_{2,22}(\Sigma_{2,12}^{\prime}-\Sigma
_{1,12}^{\prime})\}\nonumber\\[-2pt]
&&\quad\qquad{}+\frac{d}{n_1^2n_2^3}\operatorname{tr}^2(\Sigma_{1,11}\Sigma_{2,11})
\operatorname{tr}^2(\Sigma_{1,22}\Sigma_{2,22})+\frac{d}{n_2^5}\operatorname{tr}^2(\Sigma
_{2,11}^2)\operatorname{tr}^2(\Sigma_{2,22}^2)\nonumber
.
\end{eqnarray}

To evaluate the ratio of individual term in (\ref{stat2-result2}) to
$\operatorname{Var}^2(H_{n_1,n_2})$, respectively, we simply replace $\operatorname
{Var}^2(H_{n_1,n_2})$ by corresponding terms in (\ref
{stat2-variance}). Then we can show that
$\sum_{k=n_1+1}^{n_1+n_2}\mathrm{E}(C_{n,k}^4)/\operatorname
{Var}^2(H_{n_1,n_2})\to0$. Therefore, we complete the proof\vspace*{1pt} of Lemma \ref{lem5}.
\end{pf}

With two sufficient conditions given in Lemma \ref{lem4} and \ref{lem5}, we know that
\begin{eqnarray*}
\frac{H_{n_1,n_2}-\mathrm{E}(H_{n_1,n_2})}{\operatorname{Var}(H_{n_1,n_2})}
\stackrel{d}{\rightarrow}  \mathrm{N}(0,1).
\end{eqnarray*}

If we let\vspace*{1pt} $\varepsilon
_{n_1,n_2}=U_{n_1,2}+U_{n_1,3}+U_{n_2,2}+U_{n_2,3}-2W_{n_1n_1,2}-2W_{n_1n_1,3}-2W_{n_1n_1,4}$,
then $S_{n_1,n_2}=H_{n_1,n_2}+\varepsilon_{n_1,n_2}$. From $\operatorname
{Var}(\varepsilon_{n_1,n_2})=o(\sigma_{n_1,n_2}^2)$,
\begin{eqnarray*}
\operatorname{Var}\biggl(\frac{\varepsilon_{n_1,n_2}}{\sigma_{n_1,n_2}}\biggr)&=&\frac
{\operatorname{Var}(\varepsilon_{n_1,n_2})}{\sigma_{n_1,n_2}^2}
\to0.
\end{eqnarray*}

Moreover, we know $\mathrm{E}(\varepsilon_{n_1,n_2})=0$. Therefore,
$\varepsilon_{n_1,n_2}/\sigma_{n_1,n_2}\stackrel{p}{\rightarrow}  0$. From
Slutsky's Theorem, we complete the proof of Theorem \ref{teo3}.\vadjust{\goodbreak}

\subsection{\texorpdfstring{Proof of Theorem \protect\ref{teo4}}{Proof of Theorem 4}}\label{sec6.5}
Applying the trace inequality, we know that $\operatorname{tr}^2(\Sigma_{h,12}\Sigma
_{h,12}^{\prime}) \le \operatorname{tr}(\Sigma_{h,11}^2)\operatorname{tr}(\Sigma_{h,22}^2)$.
Therefore, to prove Theorem \ref{teo4}, we first consider the case where
$\operatorname{tr}^2(\Sigma_{h,12}\Sigma_{h,12}^{\prime})=O\{ \operatorname{tr}(\Sigma
_{h,11}^2)\operatorname{tr}(\Sigma_{h,22}^2)\}$. From Theorem \ref{teo2}, we can show that
$A_{n_h}^{(1)}/\operatorname{tr}(\Sigma_{h,11}^2)\stackrel{p}{\rightarrow}  1$ and
$A_{n_h}^{(2)}/\operatorname{tr}(\Sigma_{h,22}^2)\stackrel{p}{\rightarrow}  1$. Moreover, from
(\ref{var-c}), there exists a constant $d_1$ such that
\begin{eqnarray*}
\operatorname{Var}\bigl\{C_{n1 n_2}^{(i)}/\operatorname{tr}(\Sigma_{1,ii}\Sigma_{2,ii})\bigr\}
&\le&d_1\biggl(\frac{1}{n_1}+\frac{1}{n_2}\biggr) \to0,
\end{eqnarray*}
which with $\mathrm{E}(C_{n1 n_2}^{(i)})=\operatorname{tr}(\Sigma_{1,ii}\Sigma
_{2,ii})$ implies that $C_{n1 n_2}^{(i)}/\operatorname{tr}(\Sigma_{1,ii}\Sigma
_{2,ii})\stackrel{p}{\rightarrow}  1$. Similarly, using $\operatorname{tr}^2(\Sigma
_{h,12}\Sigma_{h,12}^{\prime})=O\{ \operatorname{tr}(\Sigma_{h,11}^2)\operatorname{tr}(\Sigma
_{h,22}^2)\}$, we can find a constant $d_2$ such that
\begin{eqnarray*}
&&\operatorname{Var}\{U_{n_h}/\operatorname{tr}(\Sigma_{h,12}\Sigma_{h,12}^{\prime})\}\\
&&\qquad\le
\frac{d_2}{n_h}\{1+\operatorname{tr}(\Sigma_{h,11}^2)\operatorname{tr}(\Sigma_{h,22}^2)
/\operatorname{tr}^2(\Sigma_{h,12}\Sigma_{h,12}^{\prime})\}\\
&&\qquad\to 0,
\end{eqnarray*}
which together with $\mathrm{E}(U_{n_h})=\operatorname{tr}(\Sigma_{h,12}\Sigma
_{h,12}^{\prime})$ shows that $U_{n_h}/\operatorname{tr}(\Sigma_{h,12}\Sigma
_{h,12}^{\prime})\stackrel{p}{\rightarrow}  1$ for $h=1$ or $2$. Hence, if we define
\begin{eqnarray*}
\omega^2_{0,n_1,n_2,1}&=&2\biggl(\frac{1}{n_1}+\frac{1}{n_2}\biggr)^2
\operatorname{tr}^2(\Sigma_{12}\Sigma_{12}^{\prime})\quad \mbox{and}\\
\omega^2_{0,n_1,n_2,2}&=&2\sum_{i=1}^{2}\frac{1}{n_i^2}\operatorname{tr}(\Sigma
_{i,11}^2)\operatorname{tr}(\Sigma_{i,22}^2)
+\frac{4}{n_1n_2}\operatorname{tr}(\Sigma_{1,11}\Sigma_{2,11})\operatorname{tr}(\Sigma
_{1,22}\Sigma_{2,22}),
\end{eqnarray*}
then under $H_{0 b}\dvtx \Sigma_{1,12}=\Sigma_{2,12}=\Sigma_{12}$ and
from the mapping theorem,
%
\begin{eqnarray}\label{caseI}
\hspace*{25pt}\frac{\widehat{\omega}^2_{0,n_1,n_2}}{\omega^2_{0,n_1,n_2}}&=&\frac
{\omega^2_{0,n_1,n_2,1}}{\omega^2_{0,n_1,n_2}}\frac{2({U_{n_1}}/{n_1}+{U_{n_2}}/{n_2})^2}
{\omega^2_{0,n_1,n_2,1}}\nonumber\\
&&{}+\frac{\omega^2_{0,n_1,n_2,2}}{\omega^2_{0,n_1,n_2}}
\frac{\sum_{i=1}^2\{({2}/{n_i^2})A_{n_i}^{(1)}A_{n_i}^{(2)}\}
+({4}/{(n_1n_2)})C_{n_1n_2}^{(1)}C_{n_1n_2}^{(2)}}
{\omega^2_{0,n_1,n_2,2}}\\
&\stackrel{p}{\rightarrow}&  1.\nonumber
\end{eqnarray}

Next, we consider $\operatorname{tr}^2(\Sigma_{h,12}\Sigma_{h,12}^{\prime})=o\{
\operatorname{tr}(\Sigma_{h,11}^2)\operatorname{tr}(\Sigma_{h,22}^2)\}$. If we define
\begin{eqnarray*}
\widehat{\omega}^2_{0,n_1,n_2,1}&=&2\biggl(\frac{U_{n_1}}{n_2}+\frac
{U_{n_2}}{n_1}\biggr)^2 \quad \mbox{and}\\
\widehat{\omega}^2_{0,n_1,n_2,2}&=&\sum_{i=1}^2\biggl\{\frac
{2}{n_i}A_{n_i}^{(1)}A_{n_i}^{(2)}\biggr\}
+\frac{4}{n_1n_2}C_{n_1n_2}^{(1)}C_{n_1n_2}^{(2)},
\end{eqnarray*}
then, for a given constant $\varepsilon$, we have
\begin{eqnarray*}
\mathrm{P}\biggl(\biggl|\frac{\widehat{\omega}^2_{0,n_1,n_2}}{\omega
^2_{0,n_1,n_2}}-1\biggr|>\varepsilon\biggr)\le
\mathrm{P}\biggl(\frac{\widehat{\omega}^2_{0,n_1,n_2,1}}{\omega
^2_{0,n_1,n_2}}>\varepsilon/2\biggr)+
\mathrm{P}\biggl(\biggl|\frac{\widehat{\omega}^2_{0,n_1,n_2,2}}{\omega
^2_{0,n_1,n_2}}-1\biggr|>\varepsilon/2\biggr).
\end{eqnarray*}
Thus, we only need to show $\widehat{\omega}^2_{0,n_1,n_2,1}/\omega
^2_{0,n_1,n_2}\stackrel{p}{\rightarrow}  0$ and $\widehat{\omega}^2_{0,n_1,n_2,2}/\omega^2_{0,n_1,n_2}\stackrel{p}{\rightarrow}  1$,
respectively. First of all, we know $\widehat{\omega}^2_{0,n_1,n_2,2}/\omega^2_{0,n_1,n_2}\stackrel{p}{\rightarrow}  1$ from (\ref
{caseI}). Second, there exists a constant $d_3$ such that
\begin{eqnarray*}
\mathrm{P}\biggl(\frac{\widehat{\omega}^2_{0,n_1,n_2,1}}{\omega
^2_{0,n_1,n_2}}>\frac{\varepsilon}{2}\biggr)
&\le&d_3 \Biggl[\frac{\sum_{i=1}^2\operatorname{tr}^2(\Sigma_{i,12}\Sigma
_{i,12}^{\prime})}
{\sum_{i=1}^2\operatorname{tr}(\Sigma_{i,11}^2)\operatorname{tr}(\Sigma_{i,22}^2)}\\
& &\hspace*{15pt}{}+\sum_{i=1}^{2}\biggl\{\frac{1}{n_i}
+\frac{\operatorname{tr}^2(\Sigma_{i,12}\Sigma_{i,12}^{\prime})}{n_1\operatorname{tr}(\Sigma
_{i,11}^2)\operatorname{tr}(\Sigma_{i,22}^2)} \biggr\} \Biggr],
\end{eqnarray*}
which converges to zero under $\operatorname{tr}^2(\Sigma_{i,12}\Sigma
_{i,12}^{\prime})=o\{ \operatorname{tr}(\Sigma_{i,11}^2)\operatorname{tr}(\Sigma_{i,22}^2)\}$.
Therefore, we have $\widehat{\omega}^2_{0,n_1,n_2}/\omega
^2_{0,n_1,n_2}\stackrel{p}{\rightarrow}  1$, as claimed by Theorem \ref{teo4}.

\section*{Acknowledgment}
 We thank three referees for constructive
comments and suggestions which have improved the presentation of the paper.

%

\printaddresses


\begin{thebibliography}{46}

\bibitem[\protect\citeauthoryear{Anderson}{2003}]{r1}
%
\begin{bbook}[mr]
\bauthor{\bsnm{Anderson},~\bfnm{T.~W.}\binits{T.~W.}}
(\byear{2003}).
\btitle{An Introduction to Multivariate Statistical Analysis},
\bedition{3rd} ed.
\bpublisher{Wiley}, \baddress{Hoboken, NJ}.
\bid{mr={1990662}}
\bptok{imsref}%
\end{bbook}
%
\endbibitem

\bibitem[\protect\citeauthoryear{Bai}{1993}]{r2}
%
\begin{barticle}[mr]
\bauthor{\bsnm{Bai},~\bfnm{Z.~D.}\binits{Z.~D.}}
(\byear{1993}).
\btitle{Convergence rate of expected spectral distributions of large random
matrices. {II}. {S}ample covariance matrices}.
\bjournal{Ann. Probab.}
\bvolume{21}
\bpages{649--672}.
\bid{issn={0091-1798}, mr={1217560}}
\bptok{imsref}%
\end{barticle}
%
\endbibitem

\bibitem[\protect\citeauthoryear{Bai and Saranadasa}{1996}]{r5}
%
\begin{barticle}[mr]
\bauthor{\bsnm{Bai},~\bfnm{Zhidong}\binits{Z.}} \AND
\bauthor{\bsnm{Saranadasa},~\bfnm{Hewa}\binits{H.}}
(\byear{1996}).
\btitle{Effect of high dimension: By an example of a two sample problem}.
\bjournal{Statist. Sinica}
\bvolume{6}
\bpages{311--329}.
\bid{issn={1017-0405}, mr={1399305}}
\bptok{imsref}%
\end{barticle}
%
\endbibitem

\bibitem[\protect\citeauthoryear{Bai and Silverstein}{2010}]{r6}
%
\begin{bbook}[mr]
\bauthor{\bsnm{Bai},~\bfnm{Zhidong}\binits{Z.}} \AND
\bauthor{\bsnm{Silverstein},~\bfnm{Jack~W.}\binits{J.~W.}}
(\byear{2010}).
\btitle{Spectral Analysis of Large Dimensional Random Matrices},
\bedition{2nd} ed.
\bpublisher{Springer}, \baddress{New York}.
\bid{doi={10.1007/978-1-4419-0661-8}, mr={2567175}}
\bptok{imsref}%
\end{bbook}
%
\endbibitem

\bibitem[\protect\citeauthoryear{Bai and Yin}{1993}]{r3}
%
\begin{barticle}[mr]
\bauthor{\bsnm{Bai},~\bfnm{Z.~D.}\binits{Z.~D.}} \AND
\bauthor{\bsnm{Yin},~\bfnm{Y.~Q.}\binits{Y.~Q.}}
(\byear{1993}).
\btitle{Limit of the smallest eigenvalue of a large-dimensional sample
covariance matrix}.
\bjournal{Ann. Probab.}
\bvolume{21}
\bpages{1275--1294}.
\bid{issn={0091-1798}, mr={1235416}}
\bptok{imsref}%
\end{barticle}
%
\endbibitem

\bibitem[\protect\citeauthoryear{Bai et~al.}{2009}]{r4}
%
\begin{barticle}[mr]
\bauthor{\bsnm{Bai},~\bfnm{Zhidong}\binits{Z.}},
\bauthor{\bsnm{Jiang},~\bfnm{Dandan}\binits{D.}},
\bauthor{\bsnm{Yao},~\bfnm{Jian-Feng}\binits{J.-F.}} \AND
\bauthor{\bsnm{Zheng},~\bfnm{Shurong}\binits{S.}}
(\byear{2009}).
\btitle{Corrections to {LRT} on large-dimensional covariance matrix by {RMT}}.
\bjournal{Ann. Statist.}
\bvolume{37}
\bpages{3822--3840}.
\bid{doi={10.1214/09-AOS694}, issn={0090-5364}, mr={2572444}}
\bptok{imsref}%
\end{barticle}
%
\endbibitem

\bibitem[\protect\citeauthoryear{Barry, Nobel and Wright}{2005}]{r7}
%
\begin{barticle}[pbm]
\bauthor{\bsnm{Barry},~\bfnm{William~T.}\binits{W.~T.}},
\bauthor{\bsnm{Nobel},~\bfnm{Andrew~B.}\binits{A.~B.}} \AND
\bauthor{\bsnm{Wright},~\bfnm{Fred~A.}\binits{F.~A.}}
(\byear{2005}).
\btitle{Significance analysis of functional categories in gene expression
studies: A structured permutation approach}.
\bjournal{Bioinformatics}
\bvolume{21}
\bpages{1943--1949}.
\bid{doi={10.1093/bioinformatics/bti260}, issn={1367-4803}, pii={bti260},
pmid={15647293}}
\bptok{imsref}%
\end{barticle}
%
\endbibitem

\bibitem[\protect\citeauthoryear{Benjamini and Hochberg}{1995}]{r8}
%
\begin{barticle}[mr]
\bauthor{\bsnm{Benjamini},~\bfnm{Yoav}\binits{Y.}} \AND
\bauthor{\bsnm{Hochberg},~\bfnm{Yosef}\binits{Y.}}
(\byear{1995}).
\btitle{Controlling the false discovery rate: A practical and powerful approach
to multiple testing}.
\bjournal{J. Roy. Statist. Soc. Ser. B}
\bvolume{57}
\bpages{289--300}.
\bid{issn={0035-9246}, mr={1325392}}
\bptok{imsref}%
\end{barticle}
%
\endbibitem

\bibitem[\protect\citeauthoryear{Bickel and Levina}{2008a}]{r9}
%
\begin{barticle}[mr]
\bauthor{\bsnm{Bickel},~\bfnm{Peter~J.}\binits{P.~J.}} \AND
\bauthor{\bsnm{Levina},~\bfnm{Elizaveta}\binits{E.}}
(\byear{2008}a).
\btitle{Regularized estimation of large covariance matrices}.
\bjournal{Ann. Statist.}
\bvolume{36}
\bpages{199--227}.
\bid{doi={10.1214/009053607000000758}, issn={0090-5364}, mr={2387969}}
\bptok{imsref}%
\end{barticle}
%
\endbibitem

\bibitem[\protect\citeauthoryear{Bickel and Levina}{2008b}]{r10}
%
\begin{barticle}[mr]
\bauthor{\bsnm{Bickel},~\bfnm{Peter~J.}\binits{P.~J.}} \AND
\bauthor{\bsnm{Levina},~\bfnm{Elizaveta}\binits{E.}}
(\byear{2008}b).
\btitle{Covariance regularization by thresholding}.
\bjournal{Ann. Statist.}
\bvolume{36}
\bpages{2577--2604}.
\bid{doi={10.1214/08-AOS600}, issn={0090-5364}, mr={2485008}}
\bptok{imsref}%
\end{barticle}
%
\endbibitem

\bibitem[\protect\citeauthoryear{Cai and Jiang}{2011}]{r11}
%
\begin{barticle}[mr]
\bauthor{\bsnm{Cai},~\bfnm{T.~Tony}\binits{T.~T.}} \AND
\bauthor{\bsnm{Jiang},~\bfnm{Tiefeng}\binits{T.}}
(\byear{2011}).
\btitle{Limiting laws of coherence of random matrices with
applications to
testing covariance structure and construction of compressed sensing
matrices}.
\bjournal{Ann. Statist.}
\bvolume{39}
\bpages{1496--1525}.
\bid{doi={10.1214/11-AOS879}, issn={0090-5364}, mr={2850210}}
\bptok{imsref}%
\end{barticle}
%
\endbibitem

\bibitem[\protect\citeauthoryear{Cai, Liu and Xia}{2011}]{r12}
%
\begin{bmisc}[auto:STB|2012/04/25|11:03:30]
\bauthor{\bsnm{Cai},~\bfnm{T.}\binits{T.}},
\bauthor{\bsnm{Liu},~\bfnm{W.~D.}\binits{W.~D.}} \AND
\bauthor{\bsnm{Xia},~\bfnm{Y.}\binits{Y.}}
(\byear{2011}).
\bhowpublished{Two-sample covariance matrix testing and support
recovery. Technical report, Dept. Statistics, Univ. Pennsylvania, Philadelphia, PA}.
\bptok{imsref}%
\end{bmisc}
%
\endbibitem

\bibitem[\protect\citeauthoryear{Chen and Qin}{2010}]{r13}
%
\begin{barticle}[mr]
\bauthor{\bsnm{Chen},~\bfnm{Song~Xi}\binits{S.~X.}} \AND
\bauthor{\bsnm{Qin},~\bfnm{Ying-Li}\binits{Y.-L.}}
(\byear{2010}).
\btitle{A two-sample test for high-dimensional data with applications to
gene-set testing}.
\bjournal{Ann. Statist.}
\bvolume{38}
\bpages{808--835}.
\bid{doi={10.1214/09-AOS716}, issn={0090-5364}, mr={2604697}}
\bptok{imsref}%
\end{barticle}
%
\endbibitem

\bibitem[\protect\citeauthoryear{Chen, Zhang and Zhong}{2010}]{r14}
%
\begin{barticle}[mr]
\bauthor{\bsnm{Chen},~\bfnm{Song~Xi}\binits{S.~X.}},
\bauthor{\bsnm{Zhang},~\bfnm{Li-Xin}\binits{L.-X.}} \AND
\bauthor{\bsnm{Zhong},~\bfnm{Ping-Shou}\binits{P.-S.}}
(\byear{2010}).
\btitle{Tests for high-dimensional covariance matrices}.
\bjournal{J. Amer. Statist. Assoc.}
\bvolume{105}
\bpages{810--819}.
\bid{doi={10.1198/jasa.2010.tm09560}, issn={0162-1459}, mr={2724863}}
\bptok{imsref}%
\end{barticle}
%
\endbibitem

\bibitem[\protect\citeauthoryear{Chiaretti et~al.}{2004}]{r15}
%
\begin{barticle}[auto:STB|2012/04/25|11:03:30]
\bauthor{\bsnm{Chiaretti},~\bfnm{S.}\binits{S.}},
\bauthor{\bsnm{Li},~\bfnm{X.~C.}\binits{X.~C.}},
\bauthor{\bsnm{Gentleman},~\bfnm{R.}\binits{R.}},
\bauthor{\bsnm{Vitale},~\bfnm{A.}\binits{A.}},
\bauthor{\bsnm{Vignetti},~\bfnm{M.}\binits{M.}},
\bauthor{\bsnm{Mandelli},~\bfnm{F.}\binits{F.}},
\bauthor{\bsnm{Ritz},~\bfnm{J.}\binits{J.}} \AND
\bauthor{\bsnm{Foa},~\bfnm{R.}\binits{R.}}
(\byear{2004}).
\btitle{Gene expression profile of adult T-cell acute lymphocytic leukemia
identifies distinct subsets of patients with different response to therapy
and survival}.
\bjournal{Blood}
\bvolume{103}
\bpages{2771--2778}.
\bptok{imsref}%
\end{barticle}
%
\endbibitem

\bibitem[\protect\citeauthoryear{Donoho and Jin}{2004}]{r16}
%
\begin{barticle}[mr]
\bauthor{\bsnm{Donoho},~\bfnm{David}\binits{D.}} \AND
\bauthor{\bsnm{Jin},~\bfnm{Jiashun}\binits{J.}}
(\byear{2004}).
\btitle{Higher criticism for detecting sparse heterogeneous mixtures}.
\bjournal{Ann. Statist.}
\bvolume{32}
\bpages{962--994}.
\bid{doi={10.1214/009053604000000265}, issn={0090-5364}, mr={2065195}}
\bptok{imsref}%
\end{barticle}
%
\endbibitem

\bibitem[\protect\citeauthoryear{Dudoit, Kele{\c{s}} and van~der
Laan}{2008}]{r17}
%
\begin{bincollection}[mr]
\bauthor{\bsnm{Dudoit},~\bfnm{Sandrine}\binits{S.}},
\bauthor{\bsnm{Kele{\c{s}}},~\bfnm{S{\"u}nd{\"u}z}\binits{S.}}
\AND
\bauthor{\bparticle{van~der} \bsnm{Laan},~\bfnm{Mark~J.}\binits{M.~J.}}
(\byear{2008}).
\btitle{Multiple tests of association with biological annotation metadata}.
In \bbooktitle{Probability and Statistics: Essays in Honor of {D}avid {A}.
{F}reedman}.
\bseries{Inst. Math. Stat. Collect.}
\bvolume{2}
\bpages{153--218}.
\bpublisher{IMS}, \baddress{Beachwood, OH}.
\bid{doi={10.1214/193940307000000446}, mr={2459952}}
\bptok{imsref}%
\end{bincollection}
%
\endbibitem

\bibitem[\protect\citeauthoryear{Dykstra}{1970}]{r18}
%
\begin{barticle}[auto:STB|2012/04/25|11:03:30]
\bauthor{\bsnm{Dykstra},~\bfnm{R.~L.}\binits{R.~L.}}
(\byear{1970}).
\btitle{Establishing the positive definiteness of the sample covariance
matrix}.
\bjournal{Ann. Math. Statist.}
\bvolume{41}
\bpages{2153--2154}.
\bptok{imsref}%
\end{barticle}
%
\endbibitem

\bibitem[\protect\citeauthoryear{Efron and Tibshirani}{2007}]{r19}
%
\begin{barticle}[mr]
\bauthor{\bsnm{Efron},~\bfnm{Bradley}\binits{B.}} \AND
\bauthor{\bsnm{Tibshirani},~\bfnm{Robert}\binits{R.}}
(\byear{2007}).
\btitle{On testing the significance of sets of genes}.
\bjournal{Ann. Appl. Stat.}
\bvolume{1}
\bpages{107--129}.
\bid{doi={10.1214/07-AOAS101}, issn={1932-6157}, mr={2393843}}
\bptok{imsref}%
\end{barticle}
%
\endbibitem

\bibitem[\protect\citeauthoryear{El~Karoui}{2007}]{r20}
%
\begin{barticle}[mr]
\bauthor{\bsnm{El~Karoui},~\bfnm{Noureddine}\binits{N.}}
(\byear{2007}).
\btitle{Tracy--{W}idom limit for the largest eigenvalue of a large
class of
complex sample covariance matrices}.
\bjournal{Ann. Probab.}
\bvolume{35}
\bpages{663--714}.
\bid{doi={10.1214/009117906000000917}, issn={0091-1798}, mr={2308592}}
\bptok{imsref}%
\end{barticle}
%
\endbibitem

\bibitem[\protect\citeauthoryear{Fan, Fan and Lv}{2008}]{r21}
%
\begin{barticle}[mr]
\bauthor{\bsnm{Fan},~\bfnm{Jianqing}\binits{J.}},
\bauthor{\bsnm{Fan},~\bfnm{Yingying}\binits{Y.}} \AND
\bauthor{\bsnm{Lv},~\bfnm{Jinchi}\binits{J.}}
(\byear{2008}).
\btitle{High dimensional covariance matrix estimation using a factor model}.
\bjournal{J. Econometrics}
\bvolume{147}
\bpages{186--197}.
\bid{doi={10.1016/j.jeconom.2008.09.017}, issn={0304-4076}, mr={2472991}}
\bptok{imsref}%
\end{barticle}
%
\endbibitem

\bibitem[\protect\citeauthoryear{Fan, Hall and Yao}{2007}]{r22}
%
\begin{barticle}[mr]
\bauthor{\bsnm{Fan},~\bfnm{Jianqing}\binits{J.}},
\bauthor{\bsnm{Hall},~\bfnm{Peter}\binits{P.}} \AND
\bauthor{\bsnm{Yao},~\bfnm{Qiwei}\binits{Q.}}
(\byear{2007}).
\btitle{To how many simultaneous hypothesis tests can normal,
{S}tudent's {$t$}
or bootstrap calibration be applied?}
\bjournal{J. Amer. Statist. Assoc.}
\bvolume{102}
\bpages{1282--1288}.
\bid{doi={10.1198/016214507000000969}, issn={0162-1459}, mr={2372536}}
\bptok{imsref}%
\end{barticle}
%
\endbibitem

\bibitem[\protect\citeauthoryear{Fan, Peng and Huang}{2005}]{r23}
%
\begin{barticle}[mr]
\bauthor{\bsnm{Fan},~\bfnm{Jianqing}\binits{J.}},
\bauthor{\bsnm{Peng},~\bfnm{Heng}\binits{H.}} \AND
\bauthor{\bsnm{Huang},~\bfnm{Tao}\binits{T.}}
(\byear{2005}).
\btitle{Semilinear high-dimensional model for normalization of
microarray data:
A theoretical analysis and partial consistency}.
\bjournal{J. Amer. Statist. Assoc.}
\bvolume{100}
\bpages{781--813}.
\bid{doi={10.1198/016214504000001781}, issn={0162-1459}, mr={2201010}}
\bptnote{check related}%
\bptok{imsref}%
\end{barticle}
%
\endbibitem

\bibitem[\protect\citeauthoryear{Glasser}{1961}]{r24}
%
\begin{barticle}[mr]
\bauthor{\bsnm{Glasser},~\bfnm{Gerald~J.}\binits{G.~J.}}
(\byear{1961}).
\btitle{An unbiased estimator for powers of the arithmetic mean}.
\bjournal{J. Roy. Statist. Soc. Ser. B}
\bvolume{23}
\bpages{154--159}.
\bid{issn={0035-9246}, mr={0123388}}
\bptok{imsref}%
\end{barticle}
%
\endbibitem

\bibitem[\protect\citeauthoryear{Glasser}{1962}]{r25}
%
\begin{barticle}[mr]
\bauthor{\bsnm{Glasser},~\bfnm{Gerald~J.}\binits{G.~J.}}
(\byear{1962}).
\btitle{Estimators for the product of arithmetic means}.
\bjournal{J. Roy. Statist. Soc. Ser. B}
\bvolume{24}
\bpages{180--184}.
\bid{issn={0035-9246}, mr={0137209}}
\bptok{imsref}%
\end{barticle}
%
\endbibitem

\bibitem[\protect\citeauthoryear{Hall and Jin}{2008}]{r26}
%
\begin{barticle}[mr]
\bauthor{\bsnm{Hall},~\bfnm{Peter}\binits{P.}} \AND
\bauthor{\bsnm{Jin},~\bfnm{Jiashun}\binits{J.}}
(\byear{2008}).
\btitle{Properties of higher criticism under strong dependence}.
\bjournal{Ann. Statist.}
\bvolume{36}
\bpages{381--402}.
\bid{doi={10.1214/009053607000000767}, issn={0090-5364}, mr={2387976}}
\bptok{imsref}%
\end{barticle}
%
\endbibitem

\bibitem[\protect\citeauthoryear{Huang, Wang and Zhang}{2005}]{r28}
%
\begin{barticle}[mr]
\bauthor{\bsnm{Huang},~\bfnm{Jian}\binits{J.}},
\bauthor{\bsnm{Wang},~\bfnm{Deli}\binits{D.}} \AND
\bauthor{\bsnm{Zhang},~\bfnm{Cun-Hui}\binits{C.-H.}}
(\byear{2005}).
\btitle{A two-way semilinear model for normalization and analysis of c{DNA}
microarray data}.
\bjournal{J. Amer. Statist. Assoc.}
\bvolume{100}
\bpages{814--829}.
\bid{doi={10.1198/016214504000002032}, issn={0162-1459}, mr={2201011}}
\bptok{imsref}%
\end{barticle}
%
\endbibitem

\bibitem[\protect\citeauthoryear{Huang et~al.}{2006}]{r27}
%
\begin{barticle}[mr]
\bauthor{\bsnm{Huang},~\bfnm{Jianhua~Z.}\binits{J.~Z.}},
\bauthor{\bsnm{Liu},~\bfnm{Naiping}\binits{N.}},
\bauthor{\bsnm{Pourahmadi},~\bfnm{Mohsen}\binits{M.}} \AND
\bauthor{\bsnm{Liu},~\bfnm{Linxu}\binits{L.}}
(\byear{2006}).
\btitle{Covariance matrix selection and estimation via penalised normal
likelihood}.
\bjournal{Biometrika}
\bvolume{93}
\bpages{85--98}.
\bid{doi={10.1093/biomet/93.1.85}, issn={0006-3444}, mr={2277742}}
\bptok{imsref}%
\end{barticle}
%
\endbibitem

\bibitem[\protect\citeauthoryear{Johnstone}{2001}]{r29}
%
\begin{barticle}[mr]
\bauthor{\bsnm{Johnstone},~\bfnm{Iain~M.}\binits{I.~M.}}
(\byear{2001}).
\btitle{On the distribution of the largest eigenvalue in principal components
analysis}.
\bjournal{Ann. Statist.}
\bvolume{29}
\bpages{295--327}.
\bid{doi={10.1214/aos/1009210544}, issn={0090-5364}, mr={1863961}}
\bptok{imsref}%
\end{barticle}
%
\endbibitem

\bibitem[\protect\citeauthoryear{Johnstone and Lu}{2009}]{r30}
%
\begin{barticle}[mr]
\bauthor{\bsnm{Johnstone},~\bfnm{Iain~M.}\binits{I.~M.}} \AND
\bauthor{\bsnm{Lu},~\bfnm{Arthur~Yu}\binits{A.~Y.}}
(\byear{2009}).
\btitle{On consistency and sparsity for principal components analysis
in high
dimensions}.
\bjournal{J. Amer. Statist. Assoc.}
\bvolume{104}
\bpages{682--693}.
\bid{doi={10.1198/jasa.2009.0121}, issn={0162-1459}, mr={2751448}}
\bptok{imsref}%
\end{barticle}
%
\endbibitem

\bibitem[\protect\citeauthoryear{Lam and Yao}{2011}]{r31}
%
\begin{bmisc}[auto:STB|2012/04/25|11:03:30]
\bauthor{\bsnm{Lam},~\bfnm{C.}\binits{C.}} \AND
\bauthor{\bsnm{Yao},~\bfnm{Q.}\binits{Q.}}
(\byear{2011}).
\bhowpublished{Factor modelling for high-dimensional time series: Inference for the
number of factors. \textit{Ann. Statist.} To appear}.
\bptok{imsref}%
\end{bmisc}
%
\endbibitem

\bibitem[\protect\citeauthoryear{Lam, Yao and Bathia}{2011}]{r32}
%
\begin{barticle}[mr]
\bauthor{\bsnm{Lam},~\bfnm{Clifford}\binits{C.}},
\bauthor{\bsnm{Yao},~\bfnm{Qiwei}\binits{Q.}} \AND
\bauthor{\bsnm{Bathia},~\bfnm{Neil}\binits{N.}}
(\byear{2011}).
\btitle{Estimation of latent factors for high-dimensional time series}.
\bjournal{Biometrika}
\bvolume{98}
\bpages{901--918}.
\bid{doi={10.1093/biomet/asr048}, issn={0006-3444}, mr={2860332}}
\bptok{imsref}%
\end{barticle}
%
\endbibitem

\bibitem[\protect\citeauthoryear{Lan et~al.}{2010}]{r33}
%
\begin{bmisc}[auto:STB|2012/04/25|11:03:30]
\bauthor{\bsnm{Lan},~\bfnm{W.}\binits{W.}},
\bauthor{\bsnm{Luo},~\bfnm{R.}\binits{R.}},
\bauthor{\bsnm{Tsai},~\bfnm{C.}\binits{C.}},
\bauthor{\bsnm{Wang},~\bfnm{H.}\binits{H.}} \AND
\bauthor{\bsnm{Yang},~\bfnm{Y.}\binits{Y.}}
(\byear{2010}).
\bhowpublished{Testing the diagonality of a large covariance matrix in
a regression setting. Technical report, Peking Univ., China}.
\bptok{imsref}%
\end{bmisc}
%
\endbibitem

\bibitem[\protect\citeauthoryear{Ledoit and Wolf}{2002}]{r34}
%
\begin{barticle}[mr]
\bauthor{\bsnm{Ledoit},~\bfnm{Olivier}\binits{O.}} \AND
\bauthor{\bsnm{Wolf},~\bfnm{Michael}\binits{M.}}
(\byear{2002}).
\btitle{Some hypothesis tests for the covariance matrix when the
dimension is
large compared to the sample size}.
\bjournal{Ann. Statist.}
\bvolume{30}
\bpages{1081--1102}.
\bid{doi={10.1214/aos/1031689018}, issn={0090-5364}, mr={1926169}}
\bptok{imsref}%
\end{barticle}
%
\endbibitem

\bibitem[\protect\citeauthoryear{Ledoit and Wolf}{2004}]{r35}
%
\begin{barticle}[mr]
\bauthor{\bsnm{Ledoit},~\bfnm{Olivier}\binits{O.}} \AND
\bauthor{\bsnm{Wolf},~\bfnm{Michael}\binits{M.}}
(\byear{2004}).
\btitle{A well-conditioned estimator for large-dimensional covariance
matrices}.
\bjournal{J. Multivariate Anal.}
\bvolume{88}
\bpages{365--411}.
\bid{doi={10.1016/S0047-259X(03)00096-4}, issn={0047-259X}, mr={2026339}}
\bptok{imsref}%
\end{barticle}
%
\endbibitem

\bibitem[\protect\citeauthoryear{Nettleton, Recknor and Reecy}{2008}]{r36}
%
\begin{barticle}[pbm]
\bauthor{\bsnm{Nettleton},~\bfnm{Dan}\binits{D.}},
\bauthor{\bsnm{Recknor},~\bfnm{Justin}\binits{J.}} \AND
\bauthor{\bsnm{Reecy},~\bfnm{James~M.}\binits{J.~M.}}
(\byear{2008}).
\btitle{Identification of differentially expressed gene categories in
microarray studies using nonparametric multivariate analysis}.
\bjournal{Bioinformatics}
\bvolume{24}
\bpages{192--201}.
\bid{doi={10.1093/bioinformatics/btm583}, issn={1367-4811}, pii={btm583},
pmid={18042553}}
\bptok{imsref}%
\end{barticle}
%
\endbibitem

\bibitem[\protect\citeauthoryear{Newton et~al.}{2007}]{r37}
%
\begin{barticle}[mr]
\bauthor{\bsnm{Newton},~\bfnm{Michael~A.}\binits{M.~A.}},
\bauthor{\bsnm{Quintana},~\bfnm{Fernando~A.}\binits{F.~A.}},
\bauthor{\bparticle{den} \bsnm{Boon},~\bfnm{Johan~A.}\binits{J.~A.}},
\bauthor{\bsnm{Sengupta},~\bfnm{Srikumar}\binits{S.}} \AND
\bauthor{\bsnm{Ahlquist},~\bfnm{Paul}\binits{P.}}
(\byear{2007}).
\btitle{Random-set methods identify distinct aspects of the enrichment signal
in gene-set analysis}.
\bjournal{Ann. Appl. Stat.}
\bvolume{1}
\bpages{85--106}.
\bid{doi={10.1214/07-AOAS104}, issn={1932-6157}, mr={2393842}}
\bptok{imsref}%
\end{barticle}
%
\endbibitem

\bibitem[\protect\citeauthoryear{Rothman, Levina and Zhu}{2010}]{r38}
%
\begin{barticle}[mr]
\bauthor{\bsnm{Rothman},~\bfnm{Adam~J.}\binits{A.~J.}},
\bauthor{\bsnm{Levina},~\bfnm{Elizaveta}\binits{E.}} \AND
\bauthor{\bsnm{Zhu},~\bfnm{Ji}\binits{J.}}
(\byear{2010}).
\btitle{A new approach to {C}holesky-based covariance regularization
in high
dimensions}.
\bjournal{Biometrika}
\bvolume{97}
\bpages{539--550}.
\bid{doi={10.1093/biomet/asq022}, issn={0006-3444}, mr={2672482}}
\bptok{imsref}%
\end{barticle}
%
\endbibitem

\bibitem[\protect\citeauthoryear{Schott}{2007}]{r39}
%
\begin{barticle}[mr]
\bauthor{\bsnm{Schott},~\bfnm{James~R.}\binits{J.~R.}}
(\byear{2007}).
\btitle{A test for the equality of covariance matrices when the
dimension is
large relative to the sample sizes}.
\bjournal{Comput. Statist. Data Anal.}
\bvolume{51}
\bpages{6535--6542}.
\bid{doi={10.1016/j.csda.2007.03.004}, issn={0167-9473}, mr={2408613}}
\bptok{imsref}%
\end{barticle}
%
\endbibitem


\bibitem[\protect\citeauthoryear{Shedden and Taylor}{2004}]{r41}
%
\begin{bincollection}[auto:STB|2012/04/25|11:03:30]
\bauthor{\bsnm{Shedden},~\bfnm{K.}\binits{K.}} \AND
\bauthor{\bsnm{Taylor},~\bfnm{J.}\binits{J.}}
(\byear{2004}).
\btitle{Differential correlation detects complex associations between gene
expression and clinical outcomes in lung adenocarcinomas}.
In \bbooktitle{Methods of Microarray Data Analysis}
\bpages{IV}
(\beditor{\binits{J. S.}~\bsnm{Shoemaker}}
\AND
\beditor{\binits{S. M.} \bsnm{Lin}}, eds.)
\bpages{121--131}.
\bpublisher{Springer}, \baddress{New York}.
\bptok{imsref}%
\end{bincollection}
%
\endbibitem

\bibitem[\protect\citeauthoryear{Tracy and Widom}{1996}]{r42}
%
\begin{barticle}[mr]
\bauthor{\bsnm{Tracy},~\bfnm{Craig~A.}\binits{C.~A.}} \AND
\bauthor{\bsnm{Widom},~\bfnm{Harold}\binits{H.}}
(\byear{1996}).
\btitle{On orthogonal and symplectic matrix ensembles}.
\bjournal{Comm. Math. Phys.}
\bvolume{177}
\bpages{727--754}.
\bid{issn={0010-3616}, mr={1385083}}
\bptok{imsref}%
\end{barticle}
%
\endbibitem

\bibitem[\protect\citeauthoryear{van~der Laan and Bryan}{2001}]{r43}
%
\begin{barticle}[pbm]
\bauthor{\bparticle{van~der} \bsnm{Laan},~\bfnm{M.~J.}\binits
{M.~J.}} \AND
\bauthor{\bsnm{Bryan},~\bfnm{J.}\binits{J.}}
(\byear{2001}).
\btitle{Gene expression analysis with the parametric bootstrap}.
\bjournal{Biostatistics}
\bvolume{2}
\bpages{445--461}.
\bid{doi={10.1093/biostatistics/2.4.445}, issn={1465-4644}, pii={2/4/445},
pmid={12933635}}
\bptok{imsref}%
\end{barticle}
%
\endbibitem

\bibitem[\protect\citeauthoryear{Wu and Pourahmadi}{2003}]{r44}
%
\begin{barticle}[mr]
\bauthor{\bsnm{Wu},~\bfnm{Wei~Biao}\binits{W.~B.}} \AND
\bauthor{\bsnm{Pourahmadi},~\bfnm{Mohsen}\binits{M.}}
(\byear{2003}).
\btitle{Nonparametric estimation of large covariance matrices of longitudinal
data}.
\bjournal{Biometrika}
\bvolume{90}
\bpages{831--844}.
\bid{doi={10.1093/biomet/90.4.831}, issn={0006-3444}, mr={2024760}}
\bptok{imsref}%
\end{barticle}
%
\endbibitem

\bibitem[\protect\citeauthoryear{Zhang and Huang}{2008}]{r45}
%
\begin{barticle}[mr]
\bauthor{\bsnm{Zhang},~\bfnm{Cun-Hui}\binits{C.-H.}} \AND
\bauthor{\bsnm{Huang},~\bfnm{Jian}\binits{J.}}
(\byear{2008}).
\btitle{The sparsity and bias of the {LASSO} selection in high-dimensional
linear regression}.
\bjournal{Ann. Statist.}
\bvolume{36}
\bpages{1567--1594}.
\bid{doi={10.1214/07-AOS520}, issn={0090-5364}, mr={2435448}}
\bptok{imsref}%
\end{barticle}
%
\endbibitem

\end{thebibliography}
\end{document}